\documentclass[12pt,leqno]{article}
\usepackage{amsfonts}
\usepackage{amsmath, amsthm, amsfonts, amssymb, color}

\usepackage{mathrsfs}
\usepackage{color}
\theoremstyle{definition}

\newcommand{\scr}[1]{\mathscr #1}
\definecolor{wco}{rgb}{0.5,0.2,0.3}

\numberwithin{equation}{section} \theoremstyle{remark}

\newcommand{\ua}{\uparrow}

\title{{\bf A New Type  Distribution Dependent SDE for Singular Nonlinear PDE}\footnote{Supported in
 part by  NNSFC (11831014, 11921001) and the National Key R\&D Program of China (No. 2020YFA0712900).} }
\author{{\bf  Feng-Yu Wang   }\\
\footnotesize{ Center for Applied Mathematics, Tianjin University, Tianjin 300072, China}\\
  \footnotesize{  Department of Mathematics,
Swansea University, Bay Campus, SA1 8EN, United Kingdom}\\
\footnotesize{    wangfy@tju.edu.cn}}
\begin{document}
\allowdisplaybreaks
\def\R{\mathbb R}  \def\ff{\frac} \def\ss{\sqrt} \def\B{\mathbf
B}
\def\N{\mathbb N} \def\kk{\kappa} \def\m{{\bf m}}
\def\ee{\varepsilon}\def\ddd{D^*}
\def\dd{\delta} \def\DD{\Delta} \def\vv{\varepsilon} \def\rr{\rho}
\def\<{\langle} \def\>{\rangle}
  \def\nn{\nabla} \def\pp{\partial} \def\E{\mathbb E}
\def\d{\text{\rm{d}}} \def\bb{\beta} \def\aa{\alpha} \def\D{\scr D}
  \def\si{\sigma} \def\ess{\text{\rm{ess}}}\def\s{{\bf s}}
\def\beg{\begin} \def\beq{\begin{equation}}  \def\F{\scr F}
\def\Ric{\mathcal Ric} \def\Hess{\text{\rm{Hess}}}
\def\e{\text{\rm{e}}} \def\ua{\underline a} \def\OO{\Omega}  \def\oo{\omega}
 \def\tt{\tilde}\def\[{\lfloor} \def\]{\rfloor}
\def\cut{\text{\rm{cut}}} \def\P{\mathbb P} \def\ifn{I_n(f^{\bigotimes n})}
\def\C{\scr C}      \def\aaa{\mathbf{r}}     \def\r{r}
\def\gap{\text{\rm{gap}}} \def\prr{\pi_{{\bf m},\nu}}  \def\r{\mathbf r}
\def\Z{\mathbb Z} \def\vrr{\nu} \def\ll{\lambda}
\def\L{\scr L}\def\Tt{\tt} \def\TT{\tt}\def\II{\mathbb I}
\def\i{{\rm in}}\def\Sect{{\rm Sect}}  \def\H{\mathbb H}
\def\M{\mathbb M}\def\Q{\mathbb Q} \def\texto{\text{o}} \def\LL{\Lambda}
\def\Rank{{\rm Rank}} \def\B{\scr B} \def\i{{\rm i}} \def\HR{\hat{\R}^d}
\def\to{\rightarrow} \def\gg{\gamma}
\def\EE{\scr E} \def\W{\mathbb W}
\def\A{\scr A} \def\Lip{{\rm Lip}}\def\S{\mathbb S}
\def\BB{\scr B}\def\Ent{{\rm Ent}} \def\i{{\rm i}}\def\itparallel{{\it\parallel}}
\def\g{{\mathbf g}}\def\Sect{{\mathcal Sec}}\def\T{\mathcal T}\def\BB{{\bf B}}
\def\f\ell \def\g{\mathbf g}\def\BL{{\bf L}}  \def\BG{{\mathbb G}}
\def\Bd{{D^E}} \def\BdP{D^E_\phi} \def\Bdd{{\bf \dd}} \def\Bs{{\bf s}} \def\GA{\scr A}
\def\Bg{{\bf g}}  \def\Bdd{\psi_B} \def\supp{{\rm supp}}\def\div{{\rm div}}
\def\ddiv{{\rm div}}\def\osc{{\bf osc}}\def\1{{\bf 1}}\def\BD{\mathbb D}
\def\H{{\bf H}}\def\gg{\gamma} \def\n{{\mathbf n}}\def\GG{\Gamma}\def\T{\mathbb T}
\maketitle

\begin{abstract} We propose a new type  SDE depending on the future distributions with all initial values, and establish the correspondence between this equation   and  the associated singular nonlinear PDE. 
Well-posedness and regularities are  investigated.    \end{abstract} \noindent
 \vskip 2cm

 \section{Introduction }
Probability theory  has played a very important role in the study of partial differential equations (PDEs). Starting with the work of A. N. Kolmogorov \cite{KL}, the probabilistic approach  has been   successfully applied  in deducing scaling laws of turbulence in fluid (see \cite{FF} and references within). As a funder of stochastic calculus,   K. It\^o \cite{Ito} developed stochastic differential equations (SDEs) to construct Markov processes characterizing linear Kolmogorov-Fokker-Planck equations. Since then several probability models have been introduced to solve nonlinear PDEs, which include the backward SDE initiated by Bismut \cite{Bismut}  for the linear case and    Pardoux-Peng \cite{PP} for the nonlinear case  (see \cite{CSTV} and references within), the SDE driven by $G$-Brownian motion introduced by Peng \cite{Peng0}
(see \cite{Peng} and references within), and the distribution dependent (McKean-Vlasov) SDE introduced by McKean \cite{MC} (see \cite{BR} and references within). The first two models describe   the viscosity solution, while the third one corresponds to the weak solution (in the sense of integration by parts) of nonlinear PDEs. See also  \cite{CI, Qian}  and references within    
for characterizations on the incompressible Navier-Stokes equation    using SDEs.

 In this paper, we investigate a class of new type distribution dependent SDEs and establish correspondence to nonlinear PDEs,  such that the classical existence and uniqueness  of singular nonlinear PDEs are derived.

Let $m,d \in \mathbb N$,    and let $E= \R^d$ or $\T^d:= \R^d/\mathbb Z^d.$  For fixed $T\in (0,\infty)$, we consider   the following equation  for
$u: [0,T]\times E\to \R^m$:
\beq\label{E3} \pp_t u_t= (L_t +V_t)u_t+ F_t(\cdot, u_t,\nn u_t)\cdot \nn  u_t + g_t(\cdot, u_t,\nn u_t,\nn^2 u_t),\ \ t\in [0,T], \end{equation}
where  
$$ L_t:= {\rm tr}\{a_t\nn^2\}+b_t\cdot\nn=\sum_{i,j=1}^d a_t^{ij}\pp_i\pp_j+\sum_{i=1}^d b_t^i\pp_i$$ 
is a  time dependent second order differentiable operator on $E$ with measurable coefficients
$$a: [0,T]\times E\to \R^{d\otimes d},\ \ b: [0,T]\times E\to\R^d,$$ 
$\nn u_t=(\pp_i u_t^j)_{1\le i\le d, 1\le j\le m}$ takes values in $\R^{d\otimes m},$   $\nn^2 u_t=(\pp_i\pp_j u_t^k)_{1\le i,j\le d, 1\le k\le m}$ takes values in $\R^{d\otimes d\otimes m},$ and
\beg{align*}& g: [0,T]\times E\times \R^m\times \R^{d\otimes m}\times\R^{d\otimes d\otimes m} \to \R^m,\\
&F: [0,T]\times E\times \R^m\times \R^{d\otimes m}  \to \R^d,\ \ \ \   \ V: [0,T]\times E\to \R\end{align*}
are measurable.  
Note that formally one may put  $V_tu_t$ and $F_t \cdot \nn  u_t$ into $g_t$, but  conditions we will use for these two terms  are not covered by that for $g_t$,    see Remark 2.1 below for details.  

 To characterize \eqref{E3} using stochastic differential equation (SDE), we take   $E=\R^d$ and let $g: [0,T]\times \R^d\to\R$ only depend  on the time-space variables.
Denote 
$$\si_t(x):=\ss{2a_{T-t}(x)},\ \  t\in [0,T], x\in \R^d,$$ let  $W_t$ be a $d$-dimensional Brownian motion with respect to a complete filtration probability space $(\OO, \{\F_t\}_{t\ge 0}, \P)$, and  consider the following   SDE  on $\R^d$ where differentials are in $s\in [t,T]$: 
 \beq\label{S1}\beg{split}&\d X_{t,s}^x=\si_s(X_{t,s}^x)\d W_s
   +\big\{b_{T-s}+ F_{T-s}(\cdot,\psi_{s}, \nn\psi_{s})\big\}(X_{t,s}^x) \d s,\\
    &\psi_{s}(y):= \E\bigg[u_0(X_{s,T}^y)\e^{\int_s^T V_{T-r}(X_{s,r}^y)\d r}+ \int_s^T g_{T-r}(X_{s,r}^y)\e^{\int_s^r V_{T-r'}(X_{s,r'}^y)\d r'}\d r\bigg],\ \ y\in\R^d,\\
 &\qquad  \   s\in [t,T], \  X_{t,t}^x=x\in\R^d, \ t\in [0,T].\end{split}\end{equation}
 Since $\psi_{s}$ depends on the distribution of $\{X_{s,r}^x\}_{(r,x)\in [s,T]\times E}$,
 the SDE \eqref{S1} is distribution dependent for the future and for all initial values. This is essentially different from the usual McKean-Vlasov SDEs.
 Under certain conditions allowing coefficients to be singular in $(t,x)$, we prove the well-posedness of \eqref{S1} and make application to \eqref{E3}. 
 This extends the study of \cite{22} where  $V_t=0$ and $F_t(\cdot, u_t,\nn u_t)= -(u_t\cdot \nn)u_t$ is considered. The present model also includes the nonlinear term $F_t(\cdot,u_t, \nn u_t)=\bb |\nn u_t|^2$ in  the KPZ equation.

We will state the main results of the paper in Section 2,   explain   the   idea of solving the new type SDE and present some lemmas in Section 3,
and finally  prove the main results   in Sections 4 and 5. 

\section{Main results} 
 
We first state the main result on   \eqref{S1} and its link to PDE, then present the local and global well-posedness results on \eqref{E3}. 

\subsection{Distribution dependent SDE for nonlinear PDE}

In this part, we let $E=\R^d,$ and assume that   $g_t(x,r_1,r_2,r_3)=g_t(x)$ depends only on $(t,x)$. 
Then the PDE \eqref{E3} reduces to 
\beq\label{E2} \pp_t u_t= (L_t +V_t)u_t+ F_t(x, u_t,\nn u_t)\cdot \nn  u_t+ g_t,\ \ \ t\in [0,T].  \end{equation}
 Let   $D_T:= \{(t,s): 0\le t\le s\le T\}$. We define the solution   of \eqref{S1} as follows.

\beg{defn}  A family $X:=(X_{t,s}^x)_{(t,s,x)\in D_T\times \R^d}$  of random variables on $\R^d$ is called a solution of \eqref{S1}, if   $X_{t,s}^x$ is $\F_s$-measurable,
$(\psi_{s},\nn \psi_{s})$ exists such that
 $$\E\int_t^T  \Big\{\big|\si_s|^2+ \big|b_{T-s}+ F_{T-s}(\cdot, \psi_{s}, \nn \psi_{s})\big|\Big\}(X_{t,s}^x) \d s  <\infty$$  holds  for $ (t,x)\in [0,T]\times \R^d,$  and  $\P$-a.s.
\beg{align*}& X_{t,s}^x= x+   \int_t^s\si_r(X_{t,r}^x)\d W_r
 + \int_t^s \Big\{b_{T-r} + F_{T-s}\big(\cdot,\psi_{r},\nn\psi_{r}\big)\Big\}(X_{t,r}^x)\d r,\\
 &\ \ (t,s,x)\in D_T\times\R^d.\end{align*}  \end{defn}

To prove the well-posedness of \eqref{S1}, we introduce some conditions which allow the coefficients to be singular in $(t,x)$.
As we mentioned above that the coefficients of  \eqref{S1} may be singular in $(t,x)$. To measure the singularity,   we recall some functional spaces introduced  in \cite{XXZZ}.

 For any $p,q>1$ and $0\le t<s$, we write $f\in \tt L_q^p(t,s)$ if   $f$ is a (real or vector valued) measurable function on $[t,s]\times E$ such that
 $$\|f\|_{\tt L_q^p(t,s)}:=\sup_{z\in E} \bigg(\int_t^s\|f_r1_{B(z,1)}\|_{L^p}^q\d r\bigg)^{\ff 1 q}<\infty,$$
 where $B(z,1)$ is the unit ball at $z$, and   $\|\cdot\|_{L^p} $ is the $L^p$-norm for the Lebesgue measure.
When $f_t(x)=f(x)$ does not depend on $t$, let
  $$\|f\|_{\tt L^p}:=\sup_{z\in E}\|f1_{B(z,1)}\|_{L^p}.$$ When $E=\T^d$ which is compact, we may drop $1_{B(z,1)}$ such that $\tt L^p$ and $\tt L_q^p$ reduce to the usual 
  $L^p$ and $L_q^p$ spaces.
  
 We will take $(p,q)$ from the following class:
 $$\scr K:=\Big\{(p,q): \ p,q\in (2,\infty), \      \ff {d}p +\ff 2 q<1\Big\}.$$
Let $\|\cdot\|_\infty$ denote the uniform norm, and let $\C_b^1$ be the class of (real or vector valued) weakly differentiable functions $f$ on $E$ such that 
$$\|f\|_{\C_b^1}:=\|f\|_\infty+\|\nn f\|_\infty<\infty,$$
where $\|\nn f\|_\infty$ coincides with the Lipschitz constant of $f$. 
Let 
$$k(u_0,g,V):= \e^{\int_0^T \|V_t\|_\infty\d t} \bigg(\|u_0\|_\infty+ \int_0^T \|g_t\|_\infty\d t\bigg),$$ which might be infinite.
When $k(u_0,g,V)<\infty$,  
$$\bar F_t(x):= \sup_{|r_1|\le k(u_0,g,V)}\|F_t(x,r_1,\cdot)\|_\infty$$ 
  is finite provided $F_t(x,r_1,r_2)$ is bounded in $r_2$ but locally bounded in $r_1$. 
We assume that $V, u_0, a,b, F$ and $g$ satisfy the following  conditions.

\beg{enumerate} \item[$(H_{V,u_0})$]   
 $u_0\in \C_b^1,$ and there exists $(p_0,q_0)\in \scr K$ such that $\|V\|_{\tt L_{q_0}^{p_0}(0,T)}<\infty.$  
  \item[$(H_{a,b})$]   $a$ is invertible and positive definite,  $b$ is locally bounded,    
  $$ \sup_{t\in [0,T]}\bigg\{\sup_{x\ne y} \ff{|b_t(x)-b_t(y)|}{|x-y|}+  \|a_t\|_\infty + \|a_t^{-1}\|_\infty\bigg\} <\infty,$$
$$ \lim_{\vv\to 0} \sup_{|x-y|\le \vv, t\in [0,T]} |a_t(x)-a_t(y)|=0,$$
 and there exist $l\in \mathbb N$, $ (p_i,q_i) \in\scr K$ and $1\le f_i\in \tt L_{q_i}^{p_i}(0,T), 0\le i\le l$, such that
$$ |\nn a_t(x)|\le \sum_{i=1}^l f_i (t,x),\ \ (t,x)\in [0,T]\times E. $$
  \item[$(H_{F,g}^0)$]    $\| F(\cdot,0,0) \|_{\tt L_{q_0}^{p_0}(0,T)}+\|g\|_{\tt L_{q_0}^{p_0}(0,T)}<\infty$ for   $(p_0,q_0)\in \scr K$ in $(H_{V,u_0})$,  
and there exists a constant $K>0$ such that 
 \beg{align*}   &|F_t(x,r_1,r_2)- F_t(x,\tt r_1,\tt r_2)|
  \le K\big\{1\land  (|r_1-\tt r_1|+|r_2-\tt r_2|)\big\},\\
     &\qquad (t,x)\in [0,T]\times \R^d,\ |r_1|\lor |\tt r_1| \le k(u_0,g,V),\ r_2,\tt r_2\in \R^{d\otimes m}.\end{align*}
 \end{enumerate}

Under these conditions,    we will prove the well-posedness of \eqref{S1}.
 In this case,    let
\beg{align*} &P_{t,s}f(x):= \E[f(X_{t,s}^x)],\ \ f\in \B(\R^d),\\
&P_{t,s}^V f(x):= \E\bigg[f(X_{t,s}^x)\e^{\int_t^s V_{T-r}(X_{t,r}^x)\d r}\bigg],  \ (t,s)\in D_T,\ x\in\R^d\ f\in \B(\R^d),\\
& u_t^V(f):= \int_t^T P_{t,s}^V f_s\d s,\ \ \ f\in \B([0,T]\times \R^d) \end{align*} provided the expectations and integrals exist, where $\B(\cdot)$ is the class of measurable functions 
for a measurable space.

\beg{thm}\label{T1'} Assume $(H_{V,u_0})$, $(H_{a,b})$  and   $(H_{F,g}^0)$. Then the following assertions hold.
\beg{enumerate} \item[$(1)$] The SDE \eqref{S1} has a unique solution $X:=(X_{t,s}^x)_{(t,s,x)\in D_T\times \R^d}.$
\item[$(2)$] For any $p\in [1,\infty)$ and $(t,x)\in [0,T]\times \R^d$, the derivative
$$\nn_v X_{t,s}^x:= \lim_{\vv\downarrow 0} \ff{X_{t,s}^{x+\vv v}-X_{t,s}^x}\vv,\ \  v\in\R^d, s\in [t,T]$$
exists in $L^p(\OO\to C([t,T]; \R^d), \P),$ and there exists a constant $c(p)>0$  such that
\beq\label{A*1} \sup_{(t,x)\in [0,T]\times \R^d} \E\bigg[\sup_{s\in [t,T]} |\nn_v X_{t,s}^{x}|^p\bigg]\le c(p)|v|^p,\ \ v\in \R^d.\end{equation}
\item[$(3)$] For any $0\le t<s\le T,$   $v\in\R^d$ and $f\in \B_b(\R^d)$,
\beq\label{BS} \beg{split} &\nn_v P_{t,s} f(x):=\lim_{\vv\downarrow 0} \ff{P_{t,s} f(x+\vv v)-P_{t,s} f(x)}\vv \\
&= \ff 1 {s-t}\E\bigg[f(X_{t,s}^x)\int_t^s \Big\<\si_r^{-1}(X_{t,r}^x) \nn_v X_{t,r}^x,\ \d W_r\Big\>\bigg].\end{split}\end{equation}
\item[$(4)$] For any $p\in (1,\infty]$, there exists a constant $k_p>0$ such that for any $0\le t<s\le T$ and $f\in C_b^1(\R^d)$,
\beq\label{GRD}  |\nn P_{t,s}f| \le k_p \min\Big\{(t-s)^{-\ff 1 2}\big(P_{t,s}|f|^p\big)^{\ff 1 p},\
(P_{t,s}|\nn f|^p)^{\ff 1 p}\Big\}.\end{equation}
\item[$(5)$] For any $p\in (1,\infty]$ and $V\in \tt L_{q_0}^{p_0}(0,T)$,  there exists a constant $c>0$ determined by $ (H_{a,b})$ and $\|\bar F\|_{\tt L_{q_0}^{p_0}(0,T)}+
\|V\|_{\tt L_{q_0}^{p_0}(0,T)}$, such that 
 \beq\label{GRDD}  \beg{split} &\| P_{t,s}^V f\|_{ \C_b^1} \le c \|f\|_{\C_b^1},\ \ f\in \C_b^1,\ (t,s)\in D_T,\\
&\|\nn P_{t,s}^V f\|_\infty\le c (s-t)^{-\ff 1 2}  \|f\|_\infty,\ \ 0\le t<s\le T, \ f\in \B_b(\R^d),\\
& \|u_t^V(f)\|_\infty+\|\nn^2 u^V(f)\|_{\tt L_{q_0}^{p_0}(t,T)}\le c \|f\|_{\tt L_{q_0}^{p_0}(t,T)},\ \ t\in [0,T], \ f\in \tt L_{q_0}^{p_0}(0,T).\end{split} \end{equation}
\end{enumerate} \end{thm}

The next result provides a correspondence between solutions  of \eqref{E2} and \eqref{S1}.

\beg{thm}\label{CC}  Assume $(H_{V,u_0})$, $(H_{a,b})$  and   $(H_{F,g}^0)$.
\beg{enumerate}
\item[$(1)$]  If $(u_t)_{t\in [0,T]}$ solves \eqref{E2} in the class  $\scr U(p_0,q_0)$,   then
 \beq\label{SL}\beg{split}& u_t(x)=  \E\Big[u_0(X_{T-t,T}^x)\e^{\int_{T-t}^TV_{T-r}(X_{T-t,r}^x)\d r}\Big]\\
 &\qquad +\E \int_{T-t}^T g_{T-s}(X_{T-t,s}^x)\e^{\int_{T-t}^sV_{T-r}(X_{T-t,r}^x)\d r}\d s,\ \
    (t,x)\in [0,T]\times \R^d.\end{split}\end{equation}
 \item[$(2)$] On the other hand,     $u_t$ given by $\eqref{SL}$ solves $\eqref{E2}$ such that 
$$  \|u_t\|_{\C_b^1}+\|u\|_{\tt L_{q_0}^{p_0}(t,T)}\le c\big(\|u_0\|_{\C_b^1}+\|g\|_{\tt L_{q_0}^{p_0}(0,t)}\big),\ \ t\in [0,T]$$  
holds for some constant $c>0$ determined by $(H_{a,b})$ and $\|\bar F\|_{\tt L_{q_0}^{p_0}(0,T)}+\|V\|_{\tt L_{q_0}^{p_0}(0,T)}$. 
   
\end{enumerate}
\end{thm}

\subsection{Well-posedness   of \eqref{E3}}

In this part, we let $E=\R^d$ or $\T^d$, and allow $g_t(x,u,\nn u,\nn^2 u)$ depending on $(u,\nn u,\nn^2 u)$. 
We will   prove the existence and uniqueness of \eqref{E3}   in the following sense. 

\beg{defn} We call $(u_t)_{t\in [0,T^*)}$ a solution of \eqref{E3} in the class $\scr U (p_0,q_0)$, if
\beg{enumerate}
\item[$(1)$] $T^*\in (0, T]$,
\beq\label{*0} \sup_{t\in [0,s]}  \|u_t\|_{\C_b^1} +\|\nn^2 u\|_{\tt L_{q_0}^{p_0}(0,s)}<\infty,\ \ s\in (0,T^*),\end{equation}
and for any $t \in [0,T^*),$ 
$$u_t =u_0  +\int_0^t \Big\{(L_s+V_s) u_s  +  F_s(\cdot,u_s,\nn u_s)\cdot \nn u_s + g_s(\cdot, u_s,\nn u_s,\nn^2 u_s)\Big\}\d s.$$
\item[$(2)$] When $T^*<T$,
$\limsup_{t\uparrow T^*} \|u_t\|_{\C_b^1}=\infty.$
\end{enumerate} \end{defn}

Instead of $(H_{F,g}^0)$, we will use the following weaker conditions. Let 
$$\bar g_t(x)=|g_t(x,0,0,0)|,\ \ \ (t,x)\in [0,T]\times E.$$

 \beg{enumerate} \item[$(H_{F,g})$] 
 $\|\bar F\|_{\tt L_{q_0}^{p_0}(0,T)}+ \|\bar g\|_{\tt L_{q_0}^{p_0}(0,T)}<\infty$  for $(p_0,q_0)$ in $(H_{V,u_0})$, and there exist a constant $\aa>0$ and a map $K: \mathbb N\to (0,\infty)$ such that  
 \beg{align*} &  |F_t(x,r_1,r_2)- F_t(x,\tt r_1,\tt r_2)| +   |g_t(x,r_1,r_2, r_3)- g_t(x,\tt r_1,\tt r_2,\tt r_3)|  \\
&  \le K_n(|r_1-\tt r_1|+|r_2-\tt r_2|) +\aa |r_3-\tt r_3|,\\
   &\qquad (t,x)\in [0,T]\times E,\ |r_1|\lor |r_2|\lor |\tt r_1|\lor |\tt r_2|\le n,\ r_3,\tt r_3\in \R^{d\otimes d\otimes m}.\end{align*}
  \item[$(H_{F,g}')$] There exist   constants $C,\aa >0$ and $ \theta\in (0,\ff 1 2)$ such that  
  \beg{align*}& \sup_{(t,x)\in [0,T]\times  E} |g_{t}(x, r_1,r_2, r_3)-g_t(x,0,0,0)|\\
 & \qquad\le C (\e+|r_1|+|r_2|) \{\log (\e +|r_1|+|r_2|)\}^\theta+\aa|r_3|,\\
 &\qquad  \qquad\qquad (r_1, r_2,r_3)\in\R^m\times \R^{d\otimes m}\times\R^{d\otimes d\otimes m}.
\end{align*} 
 \end{enumerate}

We have the following result.

\beg{thm}\label{T1} Assume $(H_{V,u_0})$, $(H_{a,b})$ and $(H_{F,g})$ with small enough $\aa>0$ determined by $(H_{a,b})$ and $\|V\|_{\tt L_{q_0}^{p_0}(0,T)}+ 
\|\bar F\|_{\tt L_{q_0}^{p_0}(0,T)}$.
\beg{enumerate} \item[$(1)$]    $\eqref{E3}$ has a unique solution
in $\scr U(p_0,q_0)$.
\item[$(2)$] 
 If   $(H_{F,g}')$   holds with small enough $\aa>0$ determined by $(H_{a,b})$ and $\|V\|_{\tt L_{q_0}^{p_0}(0,T)}+ 
\|\bar F\|_{\tt L_{q_0}^{p_0}(0,T)}$, then $T^*=T$ and 
\beq\label{NN} \sup_{t\in [0,T]} \|u_t\|_{\C_b^1}+\|\nn^2 u\|_{\tt L_{q_0}^{p_0}(0,T)}<\infty.\end{equation} \end{enumerate} 
\end{thm}

\paragraph{Remark 2.1.} We do not combine $Vu$ and $F\cdot\nn u$ with $g$, since   in Theorem \ref{T1} conditions on $V$ and $F$ are not covered by that on $g$.
 To see this, we assume   $(H_{a,b})$, let $F_t(x,r_1,r_2)$ and $g_t(x,r_1,r_2,r_3)=  g_t(x,r_1,r_2)$ be locally Lipschitz continuous in $(r_1,r_2)$, and  let $\|V\|_{\tt L_{q_0}^{p_0}(0,T)}<\infty$. We consider the following two situations. 
\beg{enumerate}
\item[(a)] Let for instance  $F$ and $g$ be  bounded.  Then  Theorem \ref{T1} implies   the well-posedness and non-explosion of \eqref{E3}.
However, since $|V_t|$ may be unbounded,   $(H_{F,g}')$ does not hold if we regard $V_t r_1$ as a part of $g_t(\cdot, r_1,r_2).$  
\item[(b)]  Let for instance $ \|g\|_\infty+  \int_0^T \|V_t\|_\infty<\infty$ and $F_t(\cdot,r_1,r_2)=\bb(r_1)h_t$ for some locally bounded function $\bb$ on $\R^m$ and some 
$\R^d$-valued $h\in \tt L_{q_0}^{p_0}(0,T)$.  Then $k(u_0,V,g)<\infty$, so that 
$$\|\bar F\|_{\tt L_{q_0}^{p_0}(0,T)}\le   \|h\|_{\tt L_{q_0}^{p_0}(0,T)} \sup_{|r_1|\le k(u_0,V,g)}|\bb(r_1) |  <\infty.$$   Hence,  Theorem \ref{T1} implies the the well-posedness and non-explosion of \eqref{E3}. However, since $|\bb(r_1)|$ may have arbitrary growth in $|r_1|$, $(H_{F,g}')$ does not hold if we regard $\bb_t(r_1)h_t\cdot     r_2$ 
as a part of $g_t(\cdot, r_1,r_2)$. 
 \end{enumerate}
 
 \
 
To illustrate Theorem \ref{T1}, we consider  the following  equation \eqref{RE}. When $F=0$ it reduces to the KPZ type equation 
 $$\pp_t u_t= L_t u_t +\bb \<a_t\nn u_t,\nn u_t\>+\bar V_t,$$  and 
 for $\bb=0$ it becomes the Navier-Stokes type equation
 $$\pp_t u_t= L_t u_t + F_t(\cdot, u_t)\cdot \nn u_t+\bar V_t.$$
 
 \paragraph{Example 2.1.}  Assume $(H_{a,b})$,   let $\bb\in \R$, and  let
 $$F: [0,\infty)\times E\times \R^m \to \R^d,\ \ V: [0,\infty)\times E\to\R $$ be measurable such that 
 $$\tt F_t(x,r):=\sup_{x\in E}|F_t(x,r)-F_t(x,0)|$$ is locally bounded in $(t,r)\in [0,\infty)\times \R^m$,
   and  for some $(p_0,q_0)\in \scr K$ 
 $$\|F(\cdot, 0)\|_{\tt L_{q_0}^{p_0}(0,T)}+\|V\|_{\tt L_{q_0}^{p_0}(0,T)}<\infty,\ \ \ T\in (0,T).$$  Consider the PDE
 \beq\label{RE} \pp_t u_t= L_t u_t + \bb \<a_t\nn u_t, \nn u_t\> +  F_t\big(\cdot, \varphi_\bb(u_t)\big) \cdot\nn u_t + \bar V_t,\ \ t\ge 0,\end{equation}
 where $\<a_t\nn u_t, \nn u_t\>:=(\<a_t\nn u_t^i,\nn u_t^i\>)_{1\le i\le m}$, $\bar V_t(x):= (V_t(x),\cdots, V_t(x))\in \R^m$, and 
 $$\varphi_\bb(u_t):=  (\bb^{-1} (1- \e^{-\bb u_t^i}))_{1\le i\le m}$$ which reduces to $u_t$ when $\bb=0.$ 
If either $\|V_t\|_{\infty}$ or $\|\tt F_t\|_\infty$ is locally integrable in $t\ge 0$, then for any $u_0\in \C_b^1$, \eqref{RE}  has a unique solution $u: [0,\infty)\times E\to \R^m$ satisfying
 $$\sup_{t\in [0,T]} \|\varphi_\bb(u_t)\|_{\C_b^1}+\|\nn^2 \varphi_\bb(u)\|_{\tt L_{q_0}^{p_0}(0,T)}<\infty,\ \ T\in (0,\infty).$$
 
 \beg{proof} 
 It suffices to prove the assertion up to  an arbitrarily fixed time $T\in [0,\infty)$. 
 We only consider $\bb\ne 0$, as when $\bb=0$ the desired assertion follows from Theorem \ref{T1} for $V=0$ and $g=\bar V$.  Let $v_t=\e^{-\bb u_t}:=(\e^{-\bb u_t^i})_{1\le i\le m}.$  Then \eqref{RE} up totime $T$ s equivalent to 
\beq\label{WN}  \pp_t v_t= (L_t -\bb V_t) v_t +F_t(\cdot, \bb^{-1}(1-v_t))\cdot\nn v_t,\ \ t\in [0,T],\ v_0=\e^{-\bb u_0}.\end{equation}
If  either  $\|V_t\|_\infty$ or $\|\tt F_t\|_\infty$ is locally integrable in $t\ge 0$,
 we have $\|\bar F\|_{\tt L_{q_0}^{p_0}(0,T)}<\infty$ for 
 $$\bar F_t(x):= \sup_{|r|\le k_0} |F_t(x,\bb^{-1}(1-r))|,\ \ k_0=k(u_0,-\bb V,0):=\|u_0\|_\infty\e^{\int_0^T \|\bb V_t\|_\infty\d t}.$$   
 By Theorem \ref{T1} with $g=0$ and $-\bb V_t$ replacing $V_t$,  \eqref{WN}  has a unique solution satisfying
 $$\sup_{t\in [0,T]} \|v_t\|_{\C_b^1}+ \|\nn^2 v_t\|_{\tt L_{q_0}^{p_0}(0,T)}<\infty.$$
 Then the proof is finished. 
 \end{proof}

 
 

\section{Some lemmas}

We first explain the idea of solving \eqref{S1} by fixed point theorem.

Let  $\C_T:=C(D_T; \R^d)$ be the space of all  maps
 $$\xi=(\xi_{t,s})_{(t,s)\in D_T}: D_T\to \R^d;\ \ \xi_{t,s}\ \text{is \ continuous\ in \ }(t,s)\in D_T.$$ It is a Banach space under the uniform norm $\|\cdot\|_\infty$.

 Let $\scr P_{T}$ be the space of all  probability measures   on $\C_T$ equipped with the weak topology.
 We   define
  $\GG$ as the set of all measurable maps
$$\gg=(\gg^x)_{x\in\R^d}:   \R^d\to \scr P_T$$
such that
\beq\label{PGG} \beg{split} &\psi_{s}^\gg(x):= \int_{\C_T} u_0(\xi_{s,T}) \e^{\int_s^T V_{T-s}(\xi_{s,r})\d r}
\gg^x(\d \xi)\\
&\qquad + \int_s^T \d r \int_{\C_T} g_{T-r}(\xi_{s,r}) \e^{\int_s^r V_{T-r'}(\xi_{s,r'})\d r'} \gg^x(\d \xi)
\end{split}\end{equation} is well-defined  and the weak gradient $\nn \psi_{s}^\gg$ exists for any  $s \in [0,T].$

We will see  that for any $\gg\in \GG,$ the classical  SDE
\beq\label{GG} \d X_{t,s}^{\gg,x}= \big\{b_{T-s}+F_{T-s}(\cdot,\psi_{s}^\gg, \nn \psi_{s}^\gg)\big\}(X_{t,s}^{\gg,x})\d s + \si_s(X_{t,s}^{\gg,x})\d W_s,\ \ s\in [t,T], X_{t,t}^{\gg,x}=x \end{equation}
is well-posed. In this case,  we define  a map
\beq\label{PHI} \gg\mapsto \Phi(\gg);\ x\mapsto \Phi^x(\gg):=\L_{X^{\gg,x}},  \end{equation}
where $\L_{X^{\gg,x}}$ is the distribution of the $\C_T$-valued  random variable $X^{\gg,x}:=(X_{t,s}^{\gg,x})_{(t,s)\in D_T}.$
It is clear that   $X_{t,s}^{\gg,x}$ solves \eqref{S1} if and only if $\gg$  is a fixed point of $\Phi$, i.e. $\Phi(\gg)=\gg$.
So, for the well-posedness of \eqref{S1}, it suffices to show that $\Phi$ has a unique fixed point.

To this end, we   construct a non-empty subspace $\tt\GG$ of $\GG$, such  that $\Phi$  is a contractive map on $\tt\GG$ under a complete metric. 
According to   the definition of $\psi_t^\gg$,  the   space $\tt\GG$ should be determined by the  following functionals:
for any $(t,s,x)\in D_T\times\R^d$ and $\gg\in \GG$,   let
\beg{align*}
 &\gg_{t,s}^V(f)(x):= \int_{\C_T} f(\xi_{t,s})\e^{\int_t^sV_{T-r}(\xi_{t,r})\d r}  \gg^x(\d \xi),\ \ f\in\B(\R^d),\\
 &u_t^{\gg,V}(f)(x):= \int_t^T \gg_{t,s}^V(f_s)(x)\d s \\
 &= \int_t^T\d s \int_{\C_T} f_s(\xi_{t,s})\e^{\int_t^sV_{T-r}(\xi_{t,r})\d r}  \gg^x(\d \xi),\ \ f\in \B([t,T]\times\R^d),\end{align*} provided the integrals exist.
  The class $\tt\GG$ is defined as follows.

\beg{defn}\label{D*}
Let $\tt \GG$ be the set of  $\gg\in\GG$ satisfying
  $\gg_{t,t}^x=\dd_x$ (the Dirac measure at $x$) for $(t,x) \in [0,T]\times  \R^d$,  and 
  $$\sup_{t\in [0,T]} \|\gg_t\|_V<\infty,\ \ \|\gg_t\|_V:= \sup_{x\in \R^d} \int_{\C_T} \e^{\int_t^T |V_{T-r}(\xi_{t,r})|\d r} \gg^x(\d\xi) <\infty,$$
  $$   \sup_{t\in [0,T]} \bigg\{\sup_{\|f\|_{\C_b^1}\le 1} \|\nn\gg_{t,T}^V(f)\|_{\infty}+
   \sup_{\|f\|_{\tt L_{q_0}^{p_0}(0,T)}\le 1}  \|u_t^{\gg,V}(f)\|_{\C_b^1}\bigg\}  <\infty.$$

\end{defn}

\

To apply the fixed point theorem to $\Phi$, we introduce a complete metric on $\tt \GG$, which consists of three parts: one is the weighted   variation distance 
induced by
$$\eta_t^V(\xi):= \e^{\int_t^T |V_{T-r}(\xi_{t,r})|\d r},\ \ \xi\in C([t,T];\R^d),$$
and the other two are induced by the above functionals  $\gg_{t,T}^V$ and $u_t^{\gg,V}$. More precisely, 
for any $\gg\in\GG$ and $t\in [0,T]$, let $\gg_t=(\gg_t^x)_{x\in \R^d} $  be the marginal distribution of $\gg=(\gg^x)_{x\in\R^d}$ for 
$$\xi\mapsto\xi_t:=(\xi_{t,s})_{s\in [t,T]}\in C([t,T];\R^d),$$ i.e. $\gg_t^x:= \gg^x\circ \xi_t^{-1}, x\in\R^d.$  
Then $$\gg_t(f)(x):=\int_{C([t,T];\R^d)} f(\xi)\gg_t^x(\d \xi)=\int_{\C_T} f(\xi_{t}) \gg^x(\d\xi),\ \ x\in\R^d,\ f\in \B_b(C([t,T];\R^d)).$$
For any $\gg,\tt\gg\in \tt \GG$ and $t\in [0,T]$,  let 
\beg{align*} &\|\gg_{t}-\tt\gg_{t}\|_{V}:=\sup_{f\in \B_b(C([t,T];\R^d),\|f\|_\infty\le \eta_t^V } \|\gg_{t}   (f)- \tt\gg_{t}(f)\|_\infty,\\
&\|\gg_{t}-\tt\gg_{t}\|_{\C_b^1}:=\sup_{\|f\|_{\C_b^1}\le 1 } \|\nn\{\gg_{t,T}^V(f)-\tt\gg_{t,T}^V(f)\}\|_{\infty},\\
&\|\gg_{t}-\tt\gg_{t}\|_{p_0,q_0}:= \sup_{\|f\|_{\tt L_{q_0}^{p_0}(0,T)}\le 1}  \|u_t^{\gg,V}(f)- u_t^{\tt\gg,V}(f)\|_{\C_b^1}.\end{align*}

We will prove that for any $\ll\ge 0$, $\tt\GG$ is complete under the metric
$$\rr_\ll(\gg,\tt\gg):= \sup_{(t,s)\in D_T  } \e^{-\ll(T-t)} \big\{\|\gg_{t}-\tt\gg_{t}\|_{V}+\|\gg_{t}-\tt\gg_{t}\|_{\C_b^1}+
\|\gg_{t}-\tt\gg_{t}\|_{p_0,q_0} \big\},$$
and that $\Phi$ is contractive in $\rr_\ll$ for large $\ll>0$, so that it has a unique fixed point in $\tt\GG$. To this end, we present below some lemmas.

We first extend \cite[Theorem 2.1]{YZ} with $V=0$ to the the present setting.
For any $p,q\ge 1$, let $\tt H_{q}^{2,p}(0,T)$ be the space of
$f: [0,T]\times \R^d\to \R^d$ such that
$$\|f\|_{\tt L_{q}^p(0,T)}+\|\nn f\|_{\tt L_{q}^p(0,T)}+\|\nn^2 f\|_{\tt L_{q}^p(0,T)}<\infty.$$

 \beg{lem}\label{LN3} Assume $\|V\|_{\tt L_{q_0}^{p_0}(0,T)}<\infty$ and $(H_{a,b})$ but without the condition on $|\nn a|$. 
 Let  $b^{(0)}: [0,T]\times\R^d\to\R^d$  with $\|b^{(0)}\|_{\tt L_{q}^{p}(0,T)}<\infty$ for some  $(p,q)\in \scr K$, and let 
 $$ L_t^V := L_{T-t}+b_t^{(0)}\cdot\nn+V_{T-t},\ \ t\in [0,T].$$ 
   Then for any $\ll\ge 0$ and $f\in \tt L_{q_0}^{p_0}(0,T),$   the PDE
 \beq\label{POE} (\pp_t +   L_t^V-\ll)u_t^\ll= f_t,\ \ t\in [0,T], u_T^\ll=0\end{equation}
 has a  unique  solution in $\tt H_{q_0}^{2,p_0}(0,T)$.
 Moreover,
 for any  $0<\vv< \ff 1 2 (1-\ff d {p_0}-\ff 2 {q_0})$, there exists a constant
 $c>0$ increasing in $\|b^{(0)}\|_{\tt L_{q}^{p}(0,T)}$ such that
 \beq\label{UES}\beg{split}& (1+\ll)^{\vv}(\|u_t^\ll\|_\infty+\|\nn u_t^\ll\|_\infty)+\|\nn^2 u^\ll\|_{\tt L_{q_0}^{p_0}(t,T)} \le c \|f\|_{\tt L_{q_0}^{p_0}(t,T)},\\
 &\ \ \ll\ge 0, \ f\in \tt L_{q_0}^{p_0}(0,T),\ t\in [0,T].\end{split}\end{equation}
 \end{lem}

\beg{proof} (a) Let $L^0_t= L_{T-t}+b_t^{(0)}\cdot\nn.$ By \cite[Theorem 2.1]{YZ}, for any $\ll\ge 0$ and $f\in \tt L_{q_0}^{p_0}(0,T)$, the PDE
\beq\label{PO0} (\pp_t + L_t^0-\ll)\hat u_t^\ll = f_t,\ \ \hat u_T^\ll=0\end{equation}
has a unique solution in the $\tt H_{q_0}^{2,p_0}(0,T)$,
and   for any $0<\vv<\ff 1 2 (1-\ff d {p_0}-\ff 2 {q_0})$, there exists a constant
 $c_1>0$ increasing in $\|b^{(0)}\|_{\tt L_{q}^{p}(0,T)}$ such that
 \beq\label{UES'}\beg{split} &(1+\ll)^{\vv} \big\{\|\hat u_t^\ll\|_\infty+\|\nn \hat u_t^\ll\|_\infty\big\}+\|\nn^2 \hat u^\ll\|_{\tt L_{q_0}^{p_0}(t,T)} \\
 &\le c_1 \|f\|_{\tt L_{q_0}^{p_0}(t,T)},\ \ \ll\ge 0,\ f\in \tt L_{q_0}^{p_0}(0,T),\ t\in [0,T].\end{split}\end{equation}
 Note that in  \cite[Theorem 2.1]{YZ}  the class $\tt H_{q_0}^{2,p_0}(0,T)$ is replaced by
 $\tt W_{1,q_0}^{2,p_0,w}(0,T)$, which consists of $f\in\tt H_{q_0}^{2,p_0}(0,T)$ such that for $w(x):=(1+|x|)^{-1}$ (note that $w(x)$ therein should be $(1+|x|^2)^{-\ff p 2}$ according to
\cite[Lemma 2.3]{YZ}),
\beq\label{W} |w f|+|w\nn f|+|w\nn^2 f|+ |w \pp_t f|\in \tt L_{q_0}^{p_0}(0,T).\end{equation}
 Since $w\le 1$ and $(H_{a,b})$ implies that $a$ is bounded and $|b_t(x)|\le k(1+|x|) $ for some constant $k>0$,
   a solution of \eqref{PO0} in $\tt H_{q_0}^{p_0}(0,T)$  must satisfy \eqref{W}, so we may replace  $\tt W_{1,q_0}^{2,p_0,w}(0,T)$ by $ \tt H_{q_0}^{2,p_0}(t,T)$.

(b) To solve \eqref{POE},    for any $\eta\in \B_b([0,T]\times\R^d)$, consider the PDE
\beq\label{PLT} (\pp_t + L_t^0-\ll)  u_t^{\ll,\eta} = f_t- V_{T-t} \eta_t,\ \   u_T^{\ll,\eta}=0.\end{equation}
Let $\aa(V)=\|V\|_{\tt L_{q_0}^{p_0}(0,T)}.$ We have
 $$\|f- V_{T-\cdot} \eta\|_{\tt L_{q_0}^{p_0}(0,T)}\le \|f\|_{\tt L_{q_0}^{p_0}(0,T)} +  \aa(V) \|\eta\|_\infty.$$
 Then  (a) implies that \eqref{PLT} has a unique solution in $\tt H_{q_0}^{2, p_0}(0,T)\cap \B_b([0,T]\times \R^d)$. Thus, we obtain a map
 $$\B_b([0,T]\times \R^d)\ni \eta \mapsto u^{\ll,\eta}\in \B_b([0,T]\times \R^d).$$
 So, to see that \eqref{POE} has a unique solution in $\tt H_{q_0}^{p_0}(0,T)\cap \B_b([0,T]\times \R^d)$,
 it suffices to prove that this map has a unique fixed point.

 To this end,  for any $\eta,\tt\eta\in \B_b([0,T]\times\R^d)$, let $u^{\ll,\eta,\tt\eta}:= u^{\ll,\eta}- u^{\ll, \tt\eta}.$ Then
 $$(\pp_t + L_t^0-\ll) u_t^{\ll,\eta,\tt\eta} = V_{T-t} (\tt\eta_t-\eta_t),\ \ u_T^{\ll, \eta,\tt\eta}=0.$$
 Applying \eqref{UES'} to this equation, by (a) we find a constant $c_1>0$ such that 
$$  \|u_t^{\ll,\eta}- u_t^{\ll,\eta}\|_\infty= \| u_t^{\ll,\eta,\tt\eta}\|_\infty\\
 \le c_1 \|V_{T-t} (\tt\eta-\eta)\|_{\tt L_{q_0}^{p_0}(t,T)}.$$
 Letting
 $$\|\eta-\tt\eta\|_{\infty,n}:= \sup_{t\in [0,T]}\e^{-n(T-t)} \|\eta_t-\tt\eta_t\|_\infty,\ \ n\ge 1,$$
 this implies
 $$ \|u^{\ll,\eta}- u^{\ll,\eta}\|_{\infty,n}\le c_1 \|\eta-\tt\eta\|_{\infty,n}
 \sup_{t\in [0,T]}   |V_{T-\cdot} \e^{-n(\cdot-t)}\|_{\tt L_{q_0}^{p_0}(t,T)}.$$
 Hence, when $n$ is large enough, the map $u^{\ll, \eta}$ is contractive in the complete metric $\|\eta-\tt\eta\|_{\infty,n}$. Therefore, it has a unique fixed point as desired.

(c) It remains to prove the estimate \eqref{UES}.  Let $u^\ll$ be the unique solution of \eqref{POE} in $\tt H_{q_0}^{p_0}(0,T)$. By \eqref{UES'} we obtain
 \beq\label{AI} \beg{split}& (1+\ll)^{\vv}(\|u_t^\ll\|_\infty+\|\nn  u_t^\ll\|_\infty)+\|\nn^2  u^\ll\|_{\tt L_{q_0}^{p_0}(t,T)} \\
 &\le c_1 \|f-V_{T-\cdot} u^\ll\|_{\tt L_q^p(t,T)}\le c_1 \|f\|_{\tt L_{q_0}^{p_0}(t,T)} + c_1 \|V_{T-\cdot}u^\ll\|_{\tt L_{q_0}^{p_0}(t,T)},\ \ t\in [0,T].\end{split}\end{equation}
 Letting
 $$\|u^\ll\|_{\infty,n}:= \sup_{s\in [0,T]} \e^{-n(T-s)} \|u_s^\ll\|_\infty,\ \ n\ge 1,$$
 we derive
 $$    \e^{-n(T-t)}\|u_t^\ll\|_\infty\le c_1  \|f\|_{\tt L_{q_0}^{p_0}(t,T)}+ c_1
\|u^\ll\|_{\infty,n} \|V_{T-\cdot}\e^{-n(\cdot-t)}\|_{\tt L_{q_0}^{p_0}(t,T)}.$$
 When $n\ge 1$ is large enough, this implies
 $$ \|u^\ll\|_{\infty,n}\le c_1   \|f\|_{\tt L_{q_0}^{p_0}(t,T)}+\ff 1 2 \|u^\ll\|_{\infty,n},\ \ t\in [0,T].$$
 Therefore, there exists a constant $c_2>0$ such that
 $\|u^\ll_t\|_\infty\le c_2  \|f\|_{\tt L_{q_0}^{p_0}(t,T)}, t\in [0,T].$
 Combining this with \eqref{AI} and $\|V\|_{\tt L_{q_0}^{p_0}(0,T)}<\infty$, we finish the proof.
\end{proof}

\beg{lem} \label{N1} Assume $(H_{V,u_0})$, $(H_{a,b})$  and    $(H_{F,g}^0)$. Then for any $\gg\in \GG$, the SDE  $\eqref{GG}$ is well-posed,
  and $\Phi(\gg)$ defined in $\eqref{PHI}$ is a measurable map from $\R^d$ to $\scr P_T$.  Moreover: \beg{enumerate}
\item[$(1)$]   For any $v\in\R^d$,
$$\nn_v X_{t,s}^{\gg,x}:= \lim_{\vv\downarrow 0} \ff{X_{t,s}^{\gg, x+\vv v}-X_{t,s}^x}{\vv},\ \ s\in [t,T]$$
exists in $L^p(\OO\to C([t,T];\R^d),\P)$ for any $p\in (1,\infty)$, there exists a constant $k_p>0$  determined by $p$, $(H_{a,b})$ and $ \|\bar F\|_{\tt L_{q_0}^{p_0}(0,T)}$, such that
\beq\label{DX}\sup_{\gg\in\GG} \sup_{(t,x)\in [0,T]\times\R^d}   \E \Big[\sup_{s\in [t,T]}|\nn_v X_{t,s}^{\gg,x}|^p\Big]\le k_p |v|^p,\ \ v\in\R^d,\end{equation}
\item[$(2)$] For any $p\in (1,\infty]$ there exits a constant $k_p>0$ $k_p>0$  determined by $p$, $(H_{a,b})$ and $  \|\bar F\|_{\tt L_{q_0}^{p_0}(0,T)}$,  such that  
 $$P_{t,s}^\gg f(x):= \E[f(X_{t,s}^{\gg,x})],\ \ (t,s)\in D_T, x\in \R^d, f\in \B_b(\R^d)$$ satisfies
  \beq\label{GRG}  \beg{split}&|\nn P_{t,s}^\gg f|\le k_p \min\big\{(s-t)^{-\ff 1 2}(P_{t,s}^\gg |f|^p)^{\ff 1 p},\  (P_{t,s}^\gg |\nn f|^p)^{\ff 1 p}\big\},\\
  & \ \gg\in \GG, 0\le t<s\le T, f\in C_b^1(\R^d),\end{split}\end{equation}
and  for any $v,x\in\R^d$, any $0\le t<s\le T$, 
\beq\label{BSM} \nn_v P_{t,s}^{\gg}f(x)= \ff 1 {s-t} \E\bigg[f(X_{t,s}^{\gg,x}) \int_t^s \Big\<\si_r^{-1}   (X_{t,r}^{\gg,x}) \nn_vX_{t,r}^{\gg,x},  \d W_r\Big\> \bigg],\ \ f\in \B_b(\R^d).\end{equation}
 \end{enumerate} \end{lem}

\beg{proof}  By   $(H_{F,g}^0)$,  $b_t^\gg:= F_{T-t}(\cdot, \psi_t^\gg,\nn\psi_t^\gg)$ satisfies 
\beq\label{B0} \sup_{\gg\in\GG} \|b^\gg\|_{\tt L_{q_0}^{p_0}(0,T)}\le \|\bar F\|_{\tt L_{q_0}^{p_0}(0,T)}<\infty.\end{equation}
According to \cite[Theorem 2.1]{W21e}, $(H_{a,b})$ and \eqref{B0} imply that the SDE \eqref{GG}  is well-posed, and  \eqref{DX}, \eqref{GRG} and \eqref{BSM} hold.
Noting that the Borel $\si$-field on $\scr P_T$ is induced by the maps
$$\C_T\ni \mu\mapsto \mu_{t,s}(f):= \int_{\C_{T}} f(\xi_{t,s})\mu(\d\xi),\ \ (t,s)\in D_T, f\in C_b(\R^d),$$
by combining \eqref{DX} with the continuity of $X_{t,s}^{\gg,x}$ in $(t,s)\in D_T$,
we see that $\Phi(\gg)$ is a measurable map from $ \R^d$ to $\scr P_T$.
\end{proof}

In the following, we intend to prove $\Phi\GG\subset \tt\GG$, so that  we may apply the fixed point theorem to
the map $\Phi:\tt\GG\to\tt\GG$. To this end, we introduce the Feynman-Kac semigroup 
$$P_{t,s}^{\gg,V}f(x):=\E\Big[f(X_{t,s}^{\gg,x})\e^{\int_t^s V_{T-r}(X_{t,r}^{\gg, x})\d r}\Big],\ \ (t,s,x)\in D_T\times \R^d, f\in\B_b(\R^d),$$ where $X_{t,s}^{\gg,x}$ solves \eqref{GG} with generator 
\beq\label{LG} \tt L_t^\gg:= L_{T-t}+F_{T-t}(\cdot,\psi_t^\gg,\nn\psi_t^\gg)\cdot\nn,\ \ t\in [0,T].\end{equation}

\beg{lem} \label{N1'} Assume $(H_{V,u_0})$, $(H_{a,b})$ and     $(H_{F,g}^0)$. Then the following assertions hold.
\beg{enumerate}
\item[$(1)$] There exists a constant $c>0$ determined by $(H_{a,b})$ and $\|V\|_{\tt L_{q_0}^{p_0}(0,T)}+ \|\bar F\|_{\tt L_{q_0}^{p_0}(0,T)}$, such that for any $\gg\in \GG,$ any $f: [0,T]\times\R^d \to \R$ with $\|f\|_{\tt L_{q_0}^{p_0}(0,T)}<\infty$,
$$u_t^{\Phi(\gg),V}(f) :=\E \int_t^T f_s(X_{t,s}^x) \e^{\int_t^s V_{T-r}(X_{t,r}^x)\d r}\d s=\int_t^T P_{t,s}^{\gg,V} f_{s}\d s,\ \ t\in [0,T]$$
satisfies
\beq\label{BE} (\pp_t +\tt L_t^\gg+V_{T-t})u_t^{\Phi(\gg),V}(f)= - f_{t},\ \ t\in [0,T],\end{equation}
\beq\label{B2} \|u_t^{\Phi(\gg),f}\|_\infty+\|\nn u_t^{\Phi(\gg),f}\|_\infty+ \|\nn^2 u^{\Phi(\gg),f}\|_{\tt L_{q_0}^{p_0}(t,T)}\le c \|f\|_{\tt L_{q_0}^{p_0}(t,T)},\ \ t\in [0,T].\end{equation}
\item[$(2)$] For any  $f\in \B_b(\R^d)$ and $\gg\in\GG$, 
\beq\label{KK} \pp_t P_{t,s}^{\gg,V}f = - (\tt L_{t}^\gg +V_{T-t})P_{t,s}^{\gg,V}f,\ \ t\in [0,s], s\in (0,T],\end{equation}
and there exists a constant $c>0$ determined by $(H_{a,b})$ and $\|V\|_{\tt L_{q_0}^{p_0}(0,T)}+ \|\bar F\|_{\tt L_{q_0}^{p_0}(0,T)}$, such that
\beq\label{KK2}
\sup_{\gg\in\GG}\sup_{(t,s)\in D_T} \|P_{t,s}^{\gg,V}f\|_{\C_b^1} \le c  \|f\|_{\C_b^1}, \ \ f\in \C_b^1,\end{equation}
\beq\label{KK3} \|\nn P_{t,s}^{\gg,V}f\|_\infty\le c (s-t)^{-\ff 1 2} \|f\|_\infty,\ \ 0\le t<s\le T, f\in \B_b(\R^d).\end{equation}
\end{enumerate}\end{lem}

\beg{proof} (1) By Lemma \ref{LN3}, $(H_{a,b})$, \eqref{B0} and $\|f\|_{\tt L_{q_0}^{p_0}(0,T)}<\infty$ imply that the PDE
\beq\label{BE*} (\pp_t +\tt L_t^\gg+V_{T-t})\hat u_t^\gg= - f_t,\ \ t\in [0,T],\ \hat u_T^\gg=0\end{equation} has a unique solution
satisfying
\beq\label{B2*} \|\hat u_t^\gg\|_\infty+\|\nn \hat u_t^\gg\|_\infty+ \|\nn^2 \hat u^\gg\|_{\tt L_{q_0}^{p_0}(t,T)}\le c_0\|f\|_{\tt L_q^p(t,T)}\end{equation}
for some constant $c_0>0$ determined by $(H_{a,b})$ and $\|V\|_{\tt L_{q_0}^{p_0}(0,T)}+ \|\bar F\|_{\tt L_{q_0}^{p_0}(0,T)}$.

Next, by $(H_{a,b})$ and \eqref{B0}, the Krylov estimate in \cite[Lemma 3.2(2)]{YZ} holds,
so that as shown in the proof of \cite[Lemma 4.1(2)]{XXZZ}, we obtain Khasminskii's estimate: for some function $c: (0,\infty)\to (0,\infty)$ determined by $(H_{a,b})$ and $\|V\|_{\tt L_{q_0}^{p_0}(0,T)}+ \|\bar F\|_{\tt L_{q_0}^{p_0}(0,T)}$,
\beq\label{*DD} \sup_{\gg\in \GG}\sup_{(t,x)\in [0,T]\times \R^d}  \E \e^{\ll \int_0^t|V_{T-s}(X_{t,s}^{\gg,x})|\d s}\le c(\ll),\ \ \ll\in (0,\infty).\end{equation} 

Moreover, by \eqref{BE*} and   It\^o's  formula in \cite[Lemma 3.3]{YZ},
\beg{align*} &\d \Big[\hat u_{s}^\gg(X_{t,s}^{\gg,x})\e^{\int_t^s V_{T-r}(X_{t,r}^{\gg, x})\d r}\Big]
= \e^{\int_t^s V_{T-r}(X_{t,r}^{\gg,x})\d r} (\pp_s + \tt L_{s}^\gg +V_{T-s}) \tt u_s^\gg (X_{t,s}^{\gg,x}) \d s +\d M_s\\
&= -f_s(X_{t,s}^{\gg,x})\e^{\int_t^s V_{T-r}(X_{t,r}^{\gg, x})\d r} \d s + \d M_s,\ \ s\in [t,T]\end{align*}
holds for some martingale $M_s$. This together with $\hat u_T^\gg=0$, \eqref{B2*} and \eqref{*DD} yields 
\beg{align*} & 0=\E \Big[\hat u_T^\gg(X_{t,T}^{\gg,x})\e^{\int_t^T V_{T-r}(X_{t,r}^{\gg, x})\d r}\Big]\\
&= \hat u_t^\gg(x) - \E \int_t^T f_s(X_{t,s}^{\gg,x})\e^{\int_t^s V_{T-r}(X_{t,r}^{\gg, x})\d r}\d s=\hat u_t^\gg(x)-u_t^{\Phi(\gg),V}(f).\end{align*}
Therefore,  \eqref{BE*} and \eqref{B2*} imply  \eqref{BE} and \eqref{B2}.

(2)   By \eqref{*DD}, there exists a constant $c_0>0$ determined by $(H_{a,b})$ and $\|V\|_{\tt L_{q_0}^{p_0}(0,T)}+ \|\bar F\|_{\tt L_{q_0}^{p_0}(0,T)}$, such that
\beq\label{F0} \sup_{(t,s)\in D_T, \gg\in \GG} \|P_{t,s}^{\gg, V}f\|_\infty\le c_0\|f\|_\infty,\ \ f\in\B_b(\R^d).\end{equation}
By \cite[Theorem 2.1]{YZ}, or Lemma \ref{LN3} for $V=0$, $(H_{a,b})$ and \eqref{B0} imply that for any $s\in (0,T]$, the PDE
\beq\label{F1} (\pp_t +\tt L_t^\gg)\tt u_t= - P_{t,s}^\gg f,\ \ \tt u_s=0, t\in [0,s]\end{equation}
has a unique solution satisfying
\beq\label{F2}\sup_{t\in [0,s]}\big\{ \|\tt u_t\|_\infty+\|\nn \tt u_t\|_\infty\big\}+\|\nn^2 \tt u\|_{\tt L_{q_0}^{p_0}(0,s)}\le c_1
\|f\|_\infty\end{equation}
for some constant $c_1>0$ determined by $(H_{a,b})$ and $ \|\bar F\|_{\tt L_{q_0}^{p_0}(0,T)}$. By It\^o's formula in \cite[Lemma 3.3]{YZ},
we find a martingale $(M_r)_{r\in [t,s]}$ such that \eqref{F1} implies
$$\d \tt u_r(X_{t,r}^{\gg,x})= \d M_r- P_{r,s}^\gg f(X_{t,r}^{\gg,x})\d r,\ \ r\in [t,s],$$
so that
$$0=\E[\tt u_T(X_{t,s}^{\gg,x})]= \tt u_t(x) -\int_t^s \E[P_{r,s}^\gg f(X_{t,r}^{\gg,x})]\d r=\tt u_t(x)- (s-t)P_{t,s}^\gg f.$$
Then $P_{t,s}^\gg f= \ff{\tt u_t}{s-t}$ for $t\in [0,s)$, which together with \eqref{F1} and \eqref{F2} implies
  \beq\label{F4}(\pp_t+\tt L_t^\gg)P_{t,s}^\gg f=0,\ \ t\in [0,s],\end{equation}
and
$$\sup_{t\in [0,r)} \|P_{t,s}^\gg f\|_{\C_b^1} +  \|\nn^2 P_{\cdot,s}^\gg f\|_{\tt L_{q_0}^{p_0}(0,r)}<\infty,\ \ r\in [0,s).$$
By It\^o's formula in \cite[Lemma 3.3]{YZ}, we find a martingale $\tt M_r$ such that \eqref{F4} implies
$$\d \big\{\e^{\int_t^rV_{T-r'}(X_{t,r'}^{\gg,x})\d r'} P_{r,s}^\gg f(X_{t,r}^{\gg,x})\big\}=\d \tt M_r+
\e^{\int_t^rV_{T-r'}(X_{t,r'}^{\gg,x})\d r'}\big\{V_{T-r}P_{r,s}^\gg f\big\}(X_{t,r}^{\gg,x})\d r,\ \ r\in [t,s).$$
Thus,
\beg{align*} &P_{t,s}^{\gg, V}f(x)= \lim_{r\uparrow s} \E\big[\e^{\int_t^rV_{T-r'}(X_{t,r}^{\gg,x})\d r'} P_{r,s}^\gg f(X_{t,r}^{\gg,x} ) \big]\\
&=\lim_{r\uparrow s}  \E\bigg[ P_{t,s}^\gg f(x)   +\int_t^r \e^{\int_t^{s'} V_{T-r'}(X_{t,r'}^{\gg,x})\d r'}\big\{ V_{T-s'}P_{s',s}^\gg f\big\}(X_{t,s'}^{\gg,x})  \d s'  \bigg]  \\
&= P_{t,s}^\gg f(x)+ \int_t^s P_{t,r}^{\gg,V} \big\{V_{T-r} P_{r,s}^\gg f\big\}(x)\d r,\ \ f\in C_b(\R^d).\end{align*}
By a standard approximation argument, this implies   
\beq\label{F5} P_{t,s}^{\gg,V}f= P_{t,s}^\gg f+ u_t^{\Phi(\gg),V}(\tt f),\  \ t\in [0,s],\ f\in \B_b(\R^d),\end{equation}
where $\tt f_r:= 1_{[0,s]}(r)V_{T-r} P_{r,s}^\gg f $ satisfies
$$\|\tt f\|_{\tt L_{q_0}^{p_0}(0,s)}\le \|f\|_\infty\|V\|_{\tt L_{q_0}^{p_0}(0,s)}.$$
Combining this with  \eqref{GRG} for  $p=\infty$ and  \eqref{B2}, we find a constant $c_2>0$ determined by $(H_{a,b})$ and $\|V\|_{\tt L_{q_0}^{p_0}(0,T)}+ \|\bar F\|_{\tt L_{q_0}^{p_0}(0,T)}$,  such that
$$\sup_{(t,s)\in D_T,\gg\in\GG} \|\nn P_{t,s}^{\gg,V}f\|_\infty\le c_2 \|f\|_{\C_b^1},\ \ f\in \C_b^1, $$
$$\sup_{\gg\in\GG} \|\nn P_{t,s}^{\gg,V}f\|_\infty\le c_2 (s-t)^{-\ff 1 2},\ \ f\in \B_b(\R^d), \ 0\le t<s\le T.$$
These together with \eqref{F0} implies  \eqref{KK2} and \eqref{KK3}. Finally,  by \eqref{BE}, \eqref{F4} and \eqref{F5},
we prove \eqref{KK}.
 \end{proof}

\beg{lem}\label{N2} Assume $(H_{V,u_0})$, $(H_{a,b})$  and    $(H_{F,g}^0)$.
Then $\tt\GG$ is non-empty, $\rr_\ll$-complete and satisfies  $\Phi\tt\GG \subset \tt\GG$.  \end{lem}

\beg{proof}  (a) We first prove $\Phi\GG\subset\tt\GG$.   Let $\gg\in \GG$. By   the definition  of $\Phi$ in \eqref{PHI},   we have $\Phi_{t,t}^x(\gg)=\dd_x$.
Next,
\eqref{KK2} implies that $\{\Phi(\gg)\}_{t,s}^V(f)=P_{t,s}^{\gg,V} f$  satisfies
$$ \sup_{\gg\in \GG } \| \{\Phi(\gg)\}_{t,s}^V (f) \|_{\C_b^1}\le c   \| f\|_{\C_b^1}.$$
Moreover,   by   \eqref{B2} we obtain
$$\sup_{\gg\in \GG  }\sup_{t\in [0,T]} \|u_t^{\Phi(\gg),V}(f)\|_{\C_b^1}
\le c\|f\|_{\tt L_{q_0}^{p_0}(0,T)}.$$
Thus,    $(\Phi\tt\GG \subset ) \Phi\GG  \subset \tt\GG.$

On the other hand, let  
$\gg_X$ be the distribution of the solution $X=(X_{t,s}^x)_{(t,s)\in D_T, x\in\R^d} $ to \eqref{GG} with $F=0$.   Then the above argument implies $\gg\in\tt\GG$. So, $\tt\GG\ne \emptyset$.



(b) The proof of $\rr_\ll$-completeness is more or less standard, which is  addressed below for readers' convenience.

Let $\{\gg^{(n)}\}_{n\ge 1}\subset \tt\GG$ such that
\beq\label{MN} \lim_{n,m\to\infty} \rr_\ll(\gg^{(n)},\gg^{(m)})=0.\end{equation}
Then there exists a constant $c>0$ such that
\beq\label{MN0} \sup_{n\ge 1,t\in [0,T]}\bigg\{\|\gg_t^{(n)}\|_V+ \sup_{\|f\|_{\C_b^1}\le 1}   \|\nn(\gg^{(n)})_{t,T}^V(f)\|_{\infty}+\sup_{\|f\|_{\tt L_{q_0}^{p_0}(0,T)}\le 1}  \|u_t^{\gg^{(n)},V}(f)\|_{p_0,q_0}\bigg\}\le c.\end{equation}
Moreover,
$\{\gg^{(n)}\}_{n\ge 1}$ is a Cauchy sequence under the weighted variation distance 
$$\|\gg -\tt\gg\|_{V}:=\sup_{(t,x)\in [0,T]\times \R^d} \|\gg_t-\tt\gg_t\|_{V},$$
which is complete.   
So,   there exists a unique measurable map $\gg: \R^d\to\scr P_T$ such that
\beq\label{MN2}  \lim_{n\to\infty}\sup_{t\in [0,T]} \|\gg^{(n)}_t- \gg_t\|_{V}=0.\end{equation}In particular,
\beq\label{MN4} \gg^V_{t,T} (f)(x)= \int_{\C_T} f(\xi_{t,T})\eta^V_t(\xi)  \gg_t^x(\d \xi)= \lim_{n\to\infty} (\gg^{(n)})_{t,T}^V(f)(x),\ \ f\in \B_b(\R^d).\end{equation}
We aim to prove that $\gg\in\tt\GG$ and $\rr_\ll(\gg^{(n)},\gg)\to 0$ as $n\to\infty$.

By \eqref{MN0} and \eqref{MN2}, we obtain
\beq\label{MN3}\sup_{t\in [0,T]}   \|\gg_t\|_V = \sup_{t\in [0,T]}   \lim_{n\to\infty} \|(\gg^{(n)})_{t}\|_V <\infty.\end{equation}
By
\eqref{MN0}, we find a   constant $c_1>0$ such that 
$$ \sup_{t\in [0,T]}\sup_{\|f\|_{\C_b^1}\le 1} | (\gg^{(n)})_{t,T}^V(f)(x)- (\gg^{(n)})_{t,T}^V(f)(y)|\le c_1|x-y|,\ \ x,y\in\R^d, n\ge 1.$$
Combining this with \eqref{MN4}, we obtain
$$ \sup_{t\in [0,T]}\sup_{\|f\|_{\C_b^1}\le 1}| \gg_{t,T}^V(f)(x)- \gg_{t,T}^V(f)(y)|\le c_1|x-y|,\ \ x,y\in\R^d.$$
Hence, 
$$ \sup_{t\in [0,T]}\sup_{\|f\|_{\C_b^1}\le 1} \|\nn\gg_{t,T}^V(f)\|_{\infty}\le c_1<\infty.$$
Similarly, we can  prove
 $$\sup_{t\in [0,T]}\sup_{\|f\|_{\tt L_{q_0}^{p_0}(0,T)}\le 1} \|u_t^{\gg,V}(f) \|_{\C_b^1}<\infty.$$  Therefore, by Definition \ref{D*}, we have $\gg\in\tt \GG$.

Next,  let $\|f\|_{\C_b^1}\le 1$. By \eqref{MN} and \eqref{MN4},
we find  positive constants $\{\vv_{m,n}\}_{m,n\ge 1}$   uniformly in $f$ with $\vv_{m,n}\to 0$ as $m,n\to\infty$, such that
  for any $(t,x)\in\R^d$ and $v\in\R^d$ with $|v|=1$,   
\beg{align*}&\ff 1\vv \big|\gg_{t,T}^V (f) (x+\vv v)- \gg_{t,T}^V(f)(x)-(\gg^{(n)})_{t,T}^V(f)(x+\vv v)+(\gg^{(n)})_{t,T}^V(f)(x)\big|\\
&=\lim_{m\to\infty} \ff 1\vv \big|(\gg^{(m)})_{t,T}^V (f) (x+\vv v)- (\gg^{(m)})_{t,T}^V(f)(x)-(\gg^{(n)})_{t,T}^V(f)(x+\vv v)+ (\gg^{(n)})_{t,T}^V(f)(x)\big|\\
&\le \limsup_{m\to\infty}\ff 1 \vv \int_0^\vv \big|\nn_v \big[(\gg_{t,T}^{(m)})_{t,T}^V(f)-(\gg^{(n)})_{t,T}^V(f)\big](x+ r v)\big|\d r\\
&\le \limsup_{m\to\infty}\vv_{m,n},\ \ n\ge 1,\ \vv>0, \ t\in [0,T].\end{align*}
Letting $\vv\downarrow 0$ and taking sup in $x,v\in\R^d$ with $|v|=1$,  we derive 
$$\sup_{\|f\|_{\C_b^1}\le 1}\|\nn \{\gg_{t,T}^V (f)-   (\gg^{(n)})_{t,T}^V (f)\}\|_{\infty} \le \limsup_{m\to\infty}\vv_{m,n},\ \ n\ge 1,\ t\in [0,T],$$
so that  $\vv_{m,n}\to 0$ as $m,n\to\infty$ yields
$$\lim_{n\to\infty} \sup_{t\in [0,T]}\|\gg_{t}-\gg_{t}^{(n)}\|_{\C_b^1}=0.$$
Similarly, the same holds for $\|\gg_t-\gg_t^{(n)}\|_{p_0,p_0}$. Combining these with \eqref{MN2}, we prove $\rr_\ll(\gg^{(n)},\gg)\to 0$ as $n\to\infty.$
\end{proof}

\section{Proofs of Theorem  \ref{T1'} and Theorem \ref{CC}}

\beg{proof}[Proof of Theorem \ref{T1'}]    We intend to prove that $\Phi$ has a unique fixed point $ \gg$ in $\tt\GG,$    so that \eqref{E2} has a  unique solution  given by $X_{t,s}^x=X_{t,s}^{\gg},$ and \eqref{A*1}, \eqref{BS},  \eqref{GRD} and \eqref{GRDD}   follow  from \eqref{DX}, \eqref{BSM},  \eqref{GRG}, \eqref{B2}, \eqref{KK2} and \eqref{KK3}.  

Let $\gg,\tt\gg\in \tt\GG$.
By    $(H_{F,g}^0)$, we find a constant $K>0$ such that
\beq\label{PO} \big\|F_{T-r}(\cdot,\psi_{r}^{\tt\gg},\nn\psi_{r}^{\tt\gg}) -F_{T-r}(\cdot,\psi_{r}^\gg,\nn\psi_{r}^\gg)\big\|_\infty\le K\big\{1\land \|\psi_{r}^{\tt\gg}-\psi_{r}^{\gg}\|_{\C_b^1}\big\},
 \ \ \gg,\tt\gg\in\tt\GG. \end{equation}
By   $\|u_0\|_{\C_b^1}+\|g\|_{\tt L_{q_0}^{p_0}(0,T)}<\infty$,                                         there exists a constant $c_2>0$ such that
\beq\label{P*} \beg{split}&\|\psi_{r,T}^\gg-\psi_{r,T}^{\tt\gg}\|_{\C_b^1}
 \le c_2\sup_{s\in [r,T]}  \big\{  \|\gg_{s}-\tt\gg_{s}\|_{\C_b^1} +\|\gg_{s}-\tt\gg_{s}\|_{p_0,q_0}\big\}.\end{split}\end{equation}
 Let
 $$h_r^x:= \si_r^{-1}(X_{t,r}^{\gg,x})\big\{F_{T-r}(\cdot,\psi_{r}^{\tt\gg},\nn\psi_{r}^{\tt\gg}) -F_{T-r}(\cdot,\psi_{r}^\gg,\nn\psi_{r}^\gg)\big\}(X_{t,r}^{\gg,x}),\  \ r\in [t,T].$$
 By $(H_{a,b})$, \eqref{PO} and \eqref{P*}, we find a constant $c_3>0$ such that
 \beq\label{CB} |h_r^x|\le c_2 \sup_{s\in [r,T]} 1\land \big\{\|\gg_s-\tt\gg_s\|_{\C_b^1}+\|\gg_s-\tt\gg_s\|_{p_0,q_0}\big\}=:\zeta_r.\end{equation}
 Let
 $$R^x_t:= \e^{\int_t^T \< h_r^x,\d W_r\>-\ff 1 2\int_t^T | h_r^x|^2\d r}.$$
 By Girsanov's theorem,
 $$\tt W_r:= W_r- \int_t^r h_s^x\d s,\ \ r\in [t,T]$$
 is a Brownian motion under probability $\Q_t^x:=R_t^x\P$. Reformulating \eqref{GG} as
 $$ \d X_{t,s}^{\gg,x}= \big\{b_{T-s}+F_{T-s}(\cdot,\psi_{s}^{\tt \gg}, \nn \psi_{s}^{\tt \gg})\big\}(X_{t,s}^{\gg,x})\d s + \si_s(X_{t,s}^{\gg,x})\d \tt W_s,\ \ s\in [t,T], X_{t,t}^{\gg,x}=x,$$
 by the weak uniqueness we see that the law of $X_{t,\cdot}^x:=(X_{t,s}^{\gg,x})_{s\in [t,T]}$ under $\Q_t^x$ coincides with that of $(X_{t,s}^{\tt\gg,x})_{s\in [t,T]}$ under $\P$. So, by \eqref{*DD}, we find a constant $c_3>0$ such that 
\beg{align*} &\|\Phi_t(\gg)-\Phi_t(\tt\gg)\|_{V}= \sup_{\|f\|_\infty\le \eta_t^V}\sup_{x\in\R^d}  \big|\E[(R_t^x-1)f(X_{t,\cdot}^{\gg,x})]\big|\\
&\le \sup_{x\in\R^d} \E[\eta_t^V|R_t^x-1|]\le \big(\E |\eta_t^V|^2\big)^{\ff 1 2} \big(\E |R_t^x|^2-1\big)^{\ff 1 2}\le c_3 \big(\E |R_t^x|^2-1\big)^{\ff 1 2}.\end{align*}
Noting that \eqref{CB} implies
\beg{align*}&\E |R_t^x|^2 \le \E\big[\e^{  2\int_t^T \< h_r^x,\d W_r\>-  2\int_t^T |h_r^x|^2\d r+ \int_t^T |h_r^x|^2\d r}\big] 
\le  \e^{\int_t^T |\zeta_r|^2\d r}\\
&\le 1+ \e^{\int_t^T |\zeta_r|^2\d r}\int_t^T |\zeta_r|^2\d r\le 1  +c_2^2  \e^{Tc_2^2} \int_t^T \big\{\|\gg_s-\tt\gg_s\|_{\C_b^1}+\|\gg_s-\tt\gg_s\|_{p_0,q_0}\big\}^2\d s,\end{align*} we find a constant $c_4>0$ such that
$$\|\Phi_t(\gg)-\Phi_t(\tt\gg)\|_{V}\le c_4 \bigg(\int_t^T \big\{\|\gg_s-\tt\gg_s\|_{\C_b^1}+\|\gg_s-\tt\gg_s\|_{p_0,q_0}\big\}^2\d s\bigg)^{\ff 1 2}.$$ 
  Consequently,  
\beq\label{P1} \e^{-\ll (T-t)} \|\Phi_{t}(\gg)- \Phi_{t}(\tt\gg)\|_{V} \le
  c_4 \rr_{\ll}(\gg,\tt\gg)\bigg(\int_t^T \e^{-2\ll (r-t)}\d r\bigg)^{\ff 1 2},\ \ t\in [0,T]. \end{equation}

Next,
 by  \eqref{KK} and \eqref{KK2}, for any $f\in \C_b^1(\R^d)$,  we may apply It\^o's formula in \cite[Lemma 3.3]{YZ} to deduce
\beg{align*} &\d \Big\{\e^{\int_t^r V_{T-r'}(X_{t, r'}^{\tt\gg,x})\d r'} P_{r,s}^{\gg,V} f(X_{t,r}^{\tt\gg,x}) \Big\}= \d M_r\\
&+\e^{\int_t^r V_{T-r'}(X_{t, r'}^{\tt\gg,x})\d r'}\big\{\big[F_{T-r}(\cdot,\psi_{r,T}^{\tt\gg},\nn\psi_{r,T}^{\tt\gg}) -F_{T-r}(\cdot,\psi_{r,T}^{\gg},\nn\psi_{r,T}^{\gg})\big]\cdot \nn P_{r,s}^{\gg,V} f\big\}(X_{t,r}^{\tt\gg,x})\d r
\end{align*} for some martingale $M_r, r\in [0,s]$.
Then
\beq\label{RP0} \beg{split} &P_{t,s}^{\tt \gg,V } f(x)- P_{t,s}^{\gg,V} f(x)=\E\big[\e^{\int_t^s V_{T-r}(X_{t, r}^{\tt\gg,x})\d r}P_{s,s}^{\gg,V} f(X_{t,s}^{\tt\gg,x})\big] - \E\big[\e^{\int_t^t V_{T-r}(X_{t, r}^{\tt\gg,x})\d r}P_{t,s}^{\gg,V} f(X_{t,t}^{\tt\gg,x})\big]\\
&=\int_t^s P_{t,r}^{\tt\gg,V} \big\{\big[F_{T-r}(\cdot,\psi_{r}^{\tt\gg},\nn\psi_{r}^{\tt\gg}) -F_{T-r}(\cdot,\psi_{r}^\gg,\nn\psi_{r}^\gg)\big]\cdot \nn P_{r,s}^{\gg,V} f\big\}\d r.\end{split} \end{equation}
Combining this with \eqref{KK2}, \eqref{KK3}, \eqref{PO}   and \eqref{P*},  we find constants
    $c_5,c_6>0$ such that  
\beg{align*} &\sup_{\|f\|_{\C_b^1}\le 1}\|\nn\{\Phi(\gg)\}_{t,T}^V(f)- \nn\{\Phi (\tt\gg)\}_{t,T}^V(f)\|_{\infty}  \\
 &=\sup_{\|f\|_{\C_b^1}\le 1}\|\nn P_{t,T}^{\gg,V} f-\nn P_{t,T}^{\tt\gg,V} f\|_{\infty} 
 \le c_5\int_t^T(r-t)^{-\ff 1 2}  \|\psi_{r}^\gg-\psi_{r}^{\tt\gg}\|_{\C_b^1}  \d r \\
 & \le c_6 \int_t^T(r-t)^{-\ff 1 2}\sup_{\theta\in [r,T]} \big\{\|\gg_{r}-\tt\gg_{r}\|_{\C_b^1} +\|\gg_{r}-\tt\gg_{r}\|_{p_0,q_0}\big\}\d r.\end{align*}
 Then
 \beq\label{P2} \e^{-\ll (T-t)} \| \Phi_{t}(\gg) -  \Phi_{t}(\tt\gg) \|_{\C_b^1} \le
 c_6 \rr_\ll(\gg,\tt\gg) \int_t^s(r-t)^{-\ff 1 2}\e^{-\ll(r-t)}\d r,\ \ t\in [0,T].\end{equation}


Finally, let $\|f\|_{\tt L_{q_0}^{p_0}(0,T)}\le 1$. By \eqref{BE}, for any $r\in [0,T)$, we have
\beg{align*} &(\pp_t + \tt L_t^\gg+V_{T-t}) \big\{u_t^{\Phi(\gg),V}(f)-u_t^{\Phi(\tt\gg),V}(f)\big\}= \tt f_t,\ \ t\in [0,T], u_T^{\Phi(\gg), f}- u_T^{\Phi(\tt\gg),f}=0,\\
&\tt f_t:= \big\{F_{T-t}(\cdot,\psi_t^{\tt\gg}, \nn \psi_t^{\tt\gg})- F_{T-t}(\cdot,\psi_t^\gg, \nn \psi_t^\gg)\big\}\cdot\nn u_t^{\Phi(\tt\gg),V}(f).\end{align*}
By $\|f\|_{\tt L_{q_0}^{p_0}(0,T)}\le 1$, \eqref{B2} and \eqref{PO}, there exists a constant $c>0$ such that
$$\|\tt f_t\|_\infty\le c  \|\psi_t^\gg-\psi_t^{\tt\gg}\|_{\C_b^1},\ \ t\in [0,T].$$
So, by Lemma \ref{LN3}, there exist  constants $c_7, c_8>0$ such that
\beg{align*}  \sup_{t\in [r,T]}\ \|u_t^{\Phi(\gg),V}(f)-u_t^{\Phi(\tt\gg),V}(f)\|_{\C_b^1}
 \le c_7 \|\tt f\|_{\tt L_{q_0}^{p_0}(r,T)} 
 \le c_8 \bigg(\int_r^T \|\psi_t^\gg-\psi_t^{\tt\gg}\|_{\C_b^1}^{q_0}\d t\bigg)^{\ff 1 {q_0}}.\end{align*}
Combining this with \eqref{P*}, we find a constant $c_9>0$ such that
\beg{align*} &\e^{-\ll (T-t)}  \|\{\Phi(\gg)\}_t -\{\Phi(\tt\gg)\}_t \|_{p_0,q_0}
 \le c_9 \rr_{\ll}(\gg,\tt\gg) \bigg(\int_t^T \e^{-q_0\ll(s-t)}\d s\bigg)^{\ff 1 {q_0}},\ \ t\in [0,T].\end{align*}
By taking large enough $\ll>0$, this and \eqref{P1}-\eqref{P2} yield that
 $$\rr_\ll(\Phi(\gg),\Phi(\tt\gg))\le \ff 1 2 \rr_\ll(\gg,\tt\gg),\ \ \gg,\tt\gg\in \tt\GG.$$
 Therefore, $\Phi$ has a unique fixed point in $\tt\GG$ as desired.  

  Since for any solution  $X=(X_{t,s}^x)_{(t,s,x)\in D_T\times \R^d}$  of \eqref{E2},   Lemmas \ref{N1} and  \ref{N1'} imply   that the law $\gg_X$    of $X$ is in $\tt\GG$, so that $\gg_X$ is a fixed point of $\Phi$ in $\tt\GG$. By the uniqueness of the fixed point, as well as the well-posedness of \eqref{GG}, we prove the uniqueness of \eqref{S1}.
   \end{proof} 

\beg{proof}[Proof of Theorem \ref{CC}]  (1) Let $u$ solve \eqref{E2} for $t\in [0,T]$  and satisfy \eqref{*0}, we intend to prove that $u$ satisfies \eqref{SL}.

Note that $\psi_t=\psi_t^{\gg}$ for the unique fixed point $\gg\in\tt\GG$ of $\Phi$, \eqref{KK2} for $P_{T-t,T}^{\gg,V} u_0$ and
\eqref{B2} for $f=g_{T-\cdot}$ imply that
\beq\label{VV} \sup_{t\in [0,T]} (\| \psi_t\|_\infty+\|\nn \psi_t\|_\infty)<\infty.\end{equation}
Next, let
$$\tt f_s:= F_{T-s} (\cdot,\psi_s,\nn \psi_s)-F_{T-s}(\cdot,u_{T-s},\nn u_{T-s}).$$
By $(H_{F,g}^0)$, we find a constant $c_1>0$ such that 
\beq\label{*88} |\tt f_s|\le c_1 \|\psi_s- u_{T-s}\|_{\C_b^1},\ \ s\in [0,T]. \end{equation} 
By \eqref{S1} and It\^o's formula in   \cite[Lemma 3.3]{YZ},  we obtain
\beg{align*} &\d \big\{u_{T-s}(X_{t,s}^x)\e^{\int_t^sV_{T-r}(X_{t,r}^x)\d r}\big\}= \d M_s 
 +\e^{\int_t^sV_{T-r}(X_{t,r}^x)\d r}\big[ \tt f_s\cdot \nn u_{T-s} -g_{T-s}\big]\big(X_{t,s}^x\big)\d s,\ \ s\in [t,T]\end{align*}
for  some martingale $M_s$.
Combining this with the definition of $\psi_t$ in \eqref{S1}, we derive
\beq\label{RP1}\beg{split} &\psi_t(x)- u_{T-t}(x)\\
&= \E\bigg[u_0(X_{t,T}^x)\e^{\int_t^TV_{T-r}(X_{t,r}^x)\d r}+\int_t^T g_{T-s}(X_{t,s}^x)\e^{\int_t^sV_{T-r}(X_{t,r}^x)\d r}\d s\bigg] - u_{T-t}(x)\\
&= \E\int_t^T \e^{\int_t^sV_{T-r}(X_{t,r}^x)\d r}\big[\tt f_s \cdot \nn u_{T-s} \ \Big]\big(X_{t,s}^x\big)\d s\\
&= \int_t^T P_{t,s}^V \big[\tt f_s  \cdot \nn u_{T-s}  \big](x)\d s,\ \ t\in [0,T],\end{split}\end{equation}
where $P_{t,s}^V= P_{t,s}^{\gg,V}$ for $\gg$ being the unique fixed point of $\Phi$.
Combining this with \eqref{B2}, \eqref{*88} and  $\|\nn u\|_\infty<\infty$,   
we find a constant $c_2>0$ such that
$$\|\psi_t-u_{T-t}\|_{\C_b^1} \le c_2 \bigg(\int_t^T\|\psi_s-u_{T-s}\|_{\C_b^1}^{q_0}\d s\bigg)^{\ff 1 {q_0}},\ \ t\in [0,T].$$
Since  \eqref{*0} and \eqref{VV} imply that  $\|\psi_t-u_{T-t}\|_{\C_b^1}$ is bounded in $t\in [0,T]$,
  we prove
  $\psi_t-u_{T-t}=0$ for all $t\in [0,T]$, so that the proof is finished.
 
(2)  Let $u$ be given in \eqref{SL}, and let $P_{t,s}^V=P_{t,s}^{\gg,V}$ for
the unique fixed point $\gg\in\tt\GG$ of $\Phi$. We have
\beq\label{UT} u_t= \psi_{T-t}=P_{T-t,T}^V u_0+ \int_{T-t}^T P_{T-t,s}^Vg_{T-s}\d s,\ \ t\in [0,T].\end{equation}
Let
\beg{align*}\tt L_t = L_{T-t}+ F_{T-t}(\cdot, \psi_t,\nn\psi_t)\cdot\nn,  \ t\in [0,T].\end{align*}
 By Lemma \ref{N1'}(2),
$$u_t^{(1)}:= P_{T-t,T}^Vu_0,\ \ t\in [0,T]$$ satisfies
\beq\label{PD4} \beg{split}&\pp_t u_t^{(1)}=  (\tt L_{T-t} +V_t)u_t^{(1)},\ \ t\in [0,T],\\
&\| u^{(1)}\|_\infty+\|\nn u^{(1)}\|_\infty+\|\nn^2u^{(1)}\|_{\tt L_{q_0}^{p_0}(0,T)}<\infty.\end{split}\end{equation}

On the other hand,   by Lemma \ref{N1'}(1),
$$ u_t^{(2)}:= \int_{T-t}^T P_{T-t,s}^V g_{T-s}\d s,\ \ t\in [0,T]$$
satisfies
\beq\label{PD5} \beg{split}&\pp_t u_t^{(2)}= (\tt L_{T-t} +V_t)u_t^{(2)} + g_t,\ \ t\in [0,T],\\
&\| u^{(2)}\|_\infty+\|\nn u^{(2)}\|_\infty+\|\nn^2u^{(2)}\|_{\tt L_{q_0}^{p_0}(0,T)}<\infty,\end{split}\end{equation}
which together with \eqref{UT} and \eqref{PD4} finishes the proof.
\end{proof}

\section{Proof  of Theorem \ref{T1}   }

We first observe that by extending functions from $\T^d$ to $\R^d$ periodically, i.e. letting
$$f(x+k)=f(x),\ \ x\in\T^d,\ k\in \mathbb Z^d,$$
we extend the PDE   \eqref{E3} from $E=\T^d$ to $E=\R^d$ with assumptions $(H_{V,u_0})$, $(H_{a,b}), (H_{F,g})$ and $(H_{F,g}')$ invariant. So, the existence of the extended equation implies that of the original equation. On the other hand, if the original equation has two solutions with same  initial value, then their  periodical extensions solve
the extended equation, so that the uniqueness of the extended equation implies that of the original equation. Therefore, 
it suffices to prove for $E=\R^d$.

In the following we set $E=\R^d$,  and prove Theorem \ref{T1}. We first prove under the following stronger assumption replacing  $(H_{F,g})$, then make extension by a truncation argument. 

 \beg{enumerate} \item[$(\tt H_{F,g})$] $\|\bar F\|_{\tt L_{q_0}^{p_0}(0,T)}+ \|\bar g\|_{\tt L_{q_0}^{p_0}(0,T)}<\infty$, and there exist constants $K,\aa>0$ such that 
 \beg{align*} &  |F_t(x,r_1,r_2)- F_t(x,\tt r_1,\tt r_2)| +   |g_t(x,r_1,r_2, r_3)- g_t(x,\tt r_1,\tt r_2,\tt r_3)|  \\
&  \le K(|r_1-\tt r_1|+|r_2-\tt r_2|) +\aa |r_3-\tt r_3|,\\
   &\qquad (t,x)\in [0,T]\times \R^d,\ r_1,\tt r_1\in \R^m,\  r_2,\tt r_2\in \R^{d\otimes m},\ r_3,\tt r_3\in \R^{d\otimes d\otimes m}.\end{align*}
   \end{enumerate} 
   
   \beg{lem}\label{LW} Assume $(H_{a,b})$, $(H_{V,u_0})$ and $(\tt F_{F,g})$ with small enough $\aa>0$ determined by $(H_{a,b})$ and $\|V\|_{\tt L_{q_0}^{p_0}(0,T)}+\|\bar F\|_{\tt L_{q_0}^{p_0}(0,T)}.$  The equation $\eqref{E3}$ has a unique solution in $\scr U(p_0,q_0)$ with  $T^*=T$ such that 
   $$\sup_{t\in [0,T]} \|u_t\|_{\C_b^1} +\|\nn^2 u\|_{\tt L_{q_0}^{p_0}(0,T)}<\infty.$$\end{lem}
   
   \beg{proof} We will make use of the fixed point theorem for a map $\Psi$ induced by \eqref{E3}. Let $\H$ be the Banach space of maps 
$$h: [0,T]\times \R^d\to \R^m$$
satisfying 
$$\|h\|_\H:= \sup_{t\in [0,T]} \|h_t\|_{\C_b^1}+\|\nn^2 h\|_{\tt L_{q_0}^{p_0}(0,T)}<\infty.$$
Given $h\in\H$, consider the PDE
\beq\label{PD} \pp_t u_t^h = (L_t+V_t) u_t^h + F_t (\cdot, u_t^h,\nn u_t^h)\cdot \nn u_t^h + g_t(\cdot, h_t,\nn h_t,\nn^2 h_t),\ \ t\in [0,T].\end{equation}
By $(\tt H_{F,g})$, we have $\|\bar F\|_{\tt L_{q_0}^{p_0}(0,T)}<\infty$ and 
\beq\label{SB}\beg{split} &\|g(\cdot, h,\nn h,\nn^2 h)\|_{\tt L_{q_0}^{p_0}(0,t)}\le \|\bar g\|_{\tt L_{q_0}^{p_0}(0,t)} \\
&+ K\bigg(\int_0^t \|h_s\|_{\C_b^1}^{q_0}\d s\bigg)^{\ff 1 {q_0}}+ \aa\|\nn^2 h\|_{\tt L_{q_0}^{p_0}(0,t)}<\infty,\ \ t\in [0,T],\end{split}\end{equation} 
so that $(H_{F,g}^0)$ holds for $g(\cdot, h,\nn h,\nn^2 h)$ replacing $g$. 
According to Theorem \ref{T1'} and Theorem \ref{CC}, this together with $(H_{V, u_0})$ and  $(H_{a,b})$   implies that \eqref{PD} has a unique solution in the class $\scr U(p_0,q_0)$ up to time 
$T$, such that 
\beq\label{XX} \|u^h\|_\infty\le k(u_0,V,g),\end{equation}
and for some constant $c >0$ determined by $(H_{a,b})$ and $\|V\|_{\tt L_{q_0}^{p_0}(0,T)}+\|\bar F\|_{\tt L_{q_0}^{p_0}(0,T)}$, 
\beq\label{XX0}\beg{split}& \|u_t^h\|_{\C_b^1}+ \|\nn^2 u^h\|_{\tt L_{q_0}^{p_0}(0,t)}\le c\big(\|u_0\|_{\C_b^1}+ \|g(\cdot, h,\nn h,\nn^2 h)\|_{\tt L_{q_0}^{p_0}(0,t)}\big)\\
&\le c \big(\|u_0\|_{\C_b^1}+   \|\bar g\|_{\tt L_{q_0}^{p_0}(0,t)}\big) + cK \bigg(\int_0^t \|h_s\|_{\C_b^1}^{q_0}\d s\bigg)^{\ff 1 {q_0}}+ c\aa\|\nn^2 h\|_{\tt L_{q_0}^{p_0}(0,t)},\\
&\qquad  \ \ h\in\H, t\in [0,T].\end{split} \end{equation} 
So, it suffices to prove that the map
$$\Psi: \H\ni h\mapsto u^h\in \H$$ has a unique fixed point.

To this end, we let $\aa\in (0,\ff 1 {2c})$ where $c>0$ is determined by $(H_{a,b})$ and $\|V\|_{\tt L_{q_0}^{p_0}(0,T)}+\|\bar F\|_{\tt L_{q_0}^{p_0}(0,T)}$, and denote 
  $c_0:= c \big(\|u_0\|_{\C_b^1}+   \|\bar g\|_{\tt L_{q_0}^{p_0}(0,t)}\big)$, and let
$$\H_N:=\Big\{h\in \H:\ \|h\|_{N}:=\sup_{t\in [0,T]} \e^{-Nt} \big(\|h_t\|_{\C_b^1}+\|\nn^2 h\|_{\tt L_{q_0}^{p_0}(0,t)}\big)\le N c_0\Big\},\ \ N\ge 2.$$
We will find a constant $N_0\ge 2$ such that for any $N\ge N_0$, $\Psi:\H_N\to \H_N$ has a unique fixed point.

Firstly, we take $N_0\ge 2$ such that 
$$cK\bigg(\int_0^T \e^{-q_0N_0s}\d s\bigg)^{\ff 1 {q_0}}\le \ff 1 2.$$
Combining this with \eqref{XX0} and $c\aa\le \ff 1 2$, we obtain 
$$\|\Psi(h)\|_N=\|u^h\|_N\le c_1 +\ff 1 2 \|h\|_N\le c_0 +\ff N 2 c_0\le N c_0,\ \ h\in \H_N, N\ge N_0.$$
So, $\psi:   \H_N\to \H_N$  for $N\ge N_0$.  Thus, 
\beq\label{X*}\sup_{t\in [0,T]} \|u_t^h\|_{\C_b^1}\le   \|\Psi(h)\|_\H\le c(N):= c_0 N \e^{NT},\ \ h\in \H_N, N\ge N_0.\end{equation}  
It remains to prove  that $\Psi$ has a unique fixed point in $\H_N$ for any $N\ge N_0$.  
  
For any $h,\tt h\in \H_N,$ let 
\beg{align*} f_t:=&\,\big\{F_{T-t}(\cdot, u_{T-t}^{\tt h}, \nn u_{T-t}^{\tt h})- F_{T-t}(\cdot, u_{T-t}^{h}, \nn u_{T-t}^{h})\big\}\cdot \nn u_t^{\tt h}\\
&+ g_{T-t}(\cdot, \tt h_{T-t}, \nn  \tt h_{T-t}, \nn^2  \tt h_{T-t})- g_{T-t}(\cdot, h_{T-t}, \nn h_{T-t}, \nn^2 h_{T-t}),\ \ t\in [0,T]\end{align*}
\beg{align*} L_t^h := L_{T-t}+ F_{T-t}(\cdot, u_{T-t}^{h}, \nn u_{T-t}^{h})\cdot\nn,\ \ t\in [0,T].\end{align*}
Then 
\beq\label{WR} w_t:= u_{T-t}^h - u_{T-t}^{\tt h}=\Psi_{T-t}(h)- \Psi_{T-t}(\tt h)\end{equation}   satisfies $w\in \tt H_{q_0}^{2,p_0}(0,T)$ and solves the PDE 
$$(\pp_t + L_t^h +V_{T-t}) w_t = f_t,\ \ t\in [0,T],\ w_T=0.$$ 
By Lemma \ref{LN3}, \eqref{XX},    \eqref{X*},   $(H_{V, u_0}),$  $(H_{a,b})$ and $(\tt H_{F,g})$ imply   that for some constant  $c_1 >0$ determined by $(H_{a,b})$ and $\|V\|_{\tt L_{q_0}^{p_0}(0,T)}+ \|\bar F\|_{\tt L_{q_0}^{p_0}(0,T)}$,  
\beg{align*} & \|w_t\|_{\C_b^1}+\|\nn^2 w\|_{\tt L_{q_0}^{p_0}(t,T)} \le c_1 \|f\|_{\tt L_{q_0}^{p_0}(t,T)}\\
&\le c_1  K\bigg(\int_t^T \big\{c(N)\|w_s\|_{\C_b^1}+ \|h_s-\tt h_{T-s}\|_{\C_b^1}\big\}^{q_0}\d s\bigg)^{\ff 1 {q_0}} \\
&\quad + c_1 \aa    \|\nn^2(h_{T-\cdot}-\tt h_{T-\cdot})\|_{\tt L_{q_0}^{p_0}(t,T)},\ \ t\in [0,T].\end{align*}
Then  for any $\ll>0$, this and \eqref{WR} yield 
\beg{align*} \| \Psi(h)-\Psi(\tt h)\|_\ll  \le \vv(\ll) \big[\|h-\tt h\|_\ll  +  \|\Psi(h)-\Psi(\tt h)\|_\ll\big],\ \ h,\tt h \in\H_N \end{align*}
for $\vv(\ll):=   \max\big\{c_1\aa,\ (c(N)\lor 1)c_1K\sup_{t\in [0,T]}\big(\int_t^T \e^{-\ll q_0 (s-t)}\d s \big)^{\ff 1 {q_0}}\big\}.$ 
Taking  $\ll>0$   large enough and $\aa \in (0, \ff 1 {2c_1})$ such that 
$\vv(\ll)<\ff 1 2$, we see that    $\Psi$ is $\|\cdot\|_\ll$-contractive in $H_N$,  and hence has a unique fixed point.

\end{proof}

\beg{proof}[Proof of Theorem \ref{T1}]
(1)  For any $n\in\mathbb N$, let  $$\varphi_n(r):= r1_{\{|r|\le n\}}+ \ff{nr}{|r|}1_{\{|r|>n\}},\ \ r\in \R^m\cup\R^{d\otimes m}.$$
We take the following truncation of $F_t$:
\beg{align*}  &F_t^{(n)}(x,r_1,r_2)= F_t(x,\varphi_n(r_1), \varphi_n(r_2)),\ \ g_t^{(n)}(x,r_1,r_2, r_3)= g_t(x,\varphi_n(r_1), \varphi_n(r_2), r_3),\\
  & \ t\ge 0, x\in\R^d, (r_1,r_2, r_3)\in \R^m\times \R^{d\otimes m}\times \R^{d\otimes d\otimes d\otimes m}.\end{align*}
Then  $(H_{F,g})$ implies $(\tt H_{F,g})$ for $(F^{(n)}, g^{(n)})$ replacing $(F,g)$. So, by Lemma \ref{LW},   the equation
\beq\label{**} \beg{split} & \pp_t u_t^{(n)}= (L_t +V_t)u_t^{(n)}+ F_t^{(n)}(\cdot, u_t^{(n)}, \nn u_t^{(n)})\cdot\nn u_t^{(n)}+  g_t^{(n)}(\cdot, u_t^{(n)}, \nn u_t^{(n)}, \nn^2 u_t^{(n)}),\\
&\qquad t\in [0,T],\  u_0^{(n)}=u_0\end{split} \end{equation}  
has a unique solution in $\scr U(p_0,q_0)$ with  $T^*=T,$ and
\beq\label{C*} \sup_{t\in [0,T]} \|u_t^{(n)}\|_{\C_b^1} +\|\nn^2 u^{(n)}\|_{\tt L_{q_0}^{p_0}(0,T)} <\infty,\ \ n\ge 1.\end{equation}
Let
\beq\label{TN} T_n:=T\land \inf\{t\ge 0: \|u_t^{(n)}\|_{\C_b^1} \ge n\}.\end{equation}
Then $u_t^{(n)}$ solves \eqref{E3} up to time $T_n$. Letting $T_0=0$ and
$$T^*:=\lim_{n\to\infty} T_n,$$
   the uniqueness of \eqref{**} for every $n\ge 1$ implies  
   $$u_t^{(n)}= u_t^{(m)},\ \ t\in [0, T_{n}\land T_m],\ \ n,m\in \mathbb N,$$ so that  $\eqref{E3}$ has a unique solution
$$u_t:=\sum_{n=1}^\infty u_t^{(n)}1_{[T_{n-1}, T_n)}(t),\ \ t\in [0,T^*)$$ in the class $\scr U(p_0,q_0)$,   provided   $T^*>0$.

To prove $T^*>0$, it suffices to show that $T_n>0$ for  large    $n>0.$
For fixed $T\in (0,\infty)$, let $P_{t,s}^{(n),V}$ be the Feynman-Kac semigroup generated by
$$L_{T-t}+V_{T-t}+ F_{T-t}^{(n)}(\cdot, u_{T-t}^{(n)},\nn u_{T-t}^{(n)})\cdot\nn,\ \ t\in [0,T];$$
that is 
\beq\label{ASS}   P_{t,s}^{(n),V}f(x):= \E\big[f(X_{t,s}^{(n),x}) \e^{\int_t^s V_{T-r} (X_{t,r}^{(n),x})\d r}\big],\ \ 0\le t\le s\le T, f\in\B_b(\R^d)\end{equation} 
for $X_{t,s}^{(n),x}$ solving the SDE
\beg{align*} &\d  X_{t,s}^{(n),x}= \si_s(X_{t,s}^{(n),x})\d W_s+ \big\{b_{T-s}+F^{(n)}_{T-s}(\cdot, u_{T-s}^{(n)},\nn u_{T-s}^{(n)})\big\} (X_{t,s}^{(n),x})\d s,\\
&\qquad\qquad  s\in [t,T],\ X_{t,t}^{(n),x}=x.\end{align*} 
Since $(H_{F,g})$ implies   $\|F^{(n)}\|_{\tt L_{q_0}^{p_0}(0,T)}<\infty,$  by $(H_{a,b})$ and \cite[Theorem 2.1]{W21e}, this SDE is well-posed. 
By \eqref{**} and It\^o's formula, we find a martingale $M_s$ such that
\beg{align*} &\d \Big\{u_{T-s}^{(n)}(X_{T-t,s}^{(n),x}) \e^{\int_{T-t}^sV_{T-r}(X_{t,r}(X_{t,r}^{(n),x})\d r}\Big\}\\
&= \d M_s-g_{T-s}^{(n)}(\cdot, u_{T-s}^{(n)}, \nn u_{T-s}^{(n)},\nn^2u_{T-s}^{(n)})(X_{T-t,s}^{(n),x})\d s,\ \ s\in [T-t,T].\end{align*}  So, 
\beq\label{AS} u_t^{(n)}= P_{T-t,T}^{(n),V} u_0 +\int_{T-t}^T P_{T-t,s}^{(n),V} g_{T-s}^{(n)}(\cdot, u_{T-s}^{(n)}, \nn u_{T-s}^{(n)},\nn^2u_{T-s}^{(n)})\d s,\ \ t\in [0,T],\end{equation}  
which implies 
$\|u^{(n)}\|_\infty \le k(u_0,g,V),$  so that
$$\big|F^{(n)}_{T-s}(\cdot, u_{T-s}^{(n)},\nn u_{T-s}^{(n)})\big|\le \bar F_{T-s},\ \ s\in [0,T].$$
Combining this with $(H_{a,b})$ and $\|V\|_{\tt L_{q_0}^{p_0}(0,T)}<\infty$, 
by Lemma \ref{N1'}, we find a constant $c>0$ determined by $(H_{a,b})$ and $\|\bar F\|_{\tt L_{q_0}^{p_0}(0,T)}+ \|V\|_{\tt L_{q_0}^{p_0}(0,T)}$, such that for any $n\ge 1$ and any $t\in [0,T],$ 
\beq\label{AS1}   \|P_{T-t,T}^{(n),V} u_0\|_{\C_b^1}\le c,\ \ 
    \bigg\|\int_{T-t}^T  P_{T-t,s}^{(n),V} f_s\d s\bigg\|_{\C_b^1}\le c  \|f\|_{\tt L_{q_0}^{p_0}(T-t,T)},\ \ f\in \tt L_{q_0}^{p_0}(0,T),
  \end{equation}
 \beq\label{AS'}  \|\nn P_{t,s}^{(n),V} f\|_\infty\le c (s-t)^{-\ff 1 2}\|f\|_\infty,\ \ f\in \B_b(\R^d),\ \ s\in (t,T].\end{equation} 
 Combining   \eqref{AS} with \eqref{AS1}, and noting that $(H_{F,g})$ implies 
$$\big\|g^{(n)} (\cdot, u^{(n)},\nn u^{(n)}, \nn^2 u^{(n)})\big\|_{\tt L_{q_0}^{p_0}(0,T)}\le 2 nK_n + \|\bar g\|_{\tt L_{q_0}^{p_0}(0,T)}+\aa\|\nn^2 u^{(n)}\|_{\tt L_{q_0}^{p_0}(0,T)}<\infty,$$ we obtain 
 $$\lim_{t\downarrow 0} \|u_t^{(n)}\|_{\C_b^1}\le \lim_{t\downarrow 0}\big\{\|P_{T-t,T}^{(n),V}u_0\|_{\C_b^1} 
 + c\|g^{(n)}(\cdot, u^{(n)}, \nn u^{(n)}, \nn^2 u^{(n)})\|_{\tt L_{q_0}^{p_0}(0,t)}\big\}\le c,\ \ n\ge 1.$$
By \eqref{TN},  this implies $T_n>0$ for $n>c. $

(2) Assume $(H_{F,g}')$, we intend to prove $T^*=\infty$. For $n\ge 1$ let $ P_{t,s}^{(n),V}$ be in \eqref{ASS}. Since $u_t=u_t^{(n)}$ for $t \in [0,T_n]$,  we have
$$g_s^{(n)}(\cdot, u_s^{(n)}, \nn u_s^{(n)}, \nn^2 u_s^{(n)})= g_s (\cdot, u_s, \nn u_s, \nn^2 u_s),\ \ s\in [0,T_n].$$
Then with the integral transform $s\mapsto T-s$,  
\eqref{AS} implies
$$ u_t= P_{T-t,T}^{(n),V} u_0 +\int_{0}^t P_{T-t,T-s}^{(n),V} g_{s}(\cdot, u_{s}, \nn u_{s},\nn^2u_{s})\d s,\ \ t\in [0, T_n].$$
Combining this with \eqref{AS1}, \eqref{AS'}, we find a constant $c_1>0$ such that 
\beg{align*} &\|u_t\|_{\C_b^1}+ \|\nn^2 u\|_{\tt L_{q_0}^{p_0}(0,t)} 
 \le c\|u_0\|_{\C_b^1}  + c  \|g(\cdot, 0,0,\nn^2 u)\|_{\tt L_{q_0}^{p_0}(0,t)}\\
 &\qquad +  c\int_0^t (t-s)^{-\ff 1 2} \|g_{s}(\cdot, u_{s}, \nn u_{s},\nn^2u_{s})-g_s(\cdot,0,0,\nn^2 u_s)\|_\infty \d s.\end{align*}
 Since   $(H_{F,g}')$ implies 
 \beg{align*} &  \|g(\cdot, 0,0,\nn^2 u)\|_{\tt L_{q_0}^{p_0}(0,t)}\le \|\bar g  \|_{\tt L_{q_0}^{p_0}(0,t)}  + \aa \|\nn^2 u\|_{\tt L_{q_0}^{p_0}(0,t)},\\
 & \|g_{s}(\cdot, u_{s}, \nn u_{s},\nn^2u_{s})-g_s(\cdot,0,0,\nn^2 u_s)\|_\infty \le C  (\e+\|u_s\|_{\C_b^1}) \{\log(\e+ \|u_s\|_{\C_b^1})\}^\theta,\end{align*}
 we find a constant $c_1>0$ such that  
\beq\label{*PP} \beg{split}& \|u_t\|_{\C_b^1}+ \|\nn^2 u\|_{\tt L_{q_0}^{p_0}(0,t)}  
 \le c_1   +c\aa \|\nn^2 u\|_{\tt L_{q_0}^{p_0}(0,t)} \\
 &\quad + c_1  \int_0^t (t-s)^{-\ff 1 2} (\e+\|u_s\|_{\C_b^1}) \{\log(\e+ \|u_s\|_{\C_b^1})\}^\theta \d s,\ \ t\in [0,T_n].\end{split} \end{equation} 
Since  H\"older's inequality implies
 \beg{align*} &\int_0^t (t-s)^{-\ff 1 2}  (\e+\|u_s\|_{\C_b^1}) \{\log(\e+ \|u_s\|_{\C_b^1})\}^\theta \d s \\
&\le
\bigg(\int_0^t (t-s)^{-\ff 1 {2(1-\theta)}}\d s\bigg)^{ 1-\theta}\bigg(\int_0^t  (\e+\|u_s\|_{\C_b^1})^{\theta^{-1}}\log  (\e+ \|u_s\|_{\C_b^1})    \d s\bigg)^\theta,\end{align*} 
by  $\theta\in (0,\ff 1 2)$, when $\aa\in (0,\ff 1 c)$  we find a constant $c_2>0$ such that \eqref{*PP} implies 
\beq\label{*PP0}\beg{split} &(\e+\|u_t\|_{\C_b^1}+\|u\|_{\tt L_{q_0}^{p_0}(0,t)})^{\theta^{-1}} \\
&\le c_2+ c_2 \int_0^t (\e+\|u_s\|_{\C_b^1})^{\theta^{-1}}   \log(\e+ \|u_s\|_{\C_b^1})   \d r,\ \ t\in [0, T_n], n\ge 1.\end{split}\end{equation} 
By Gronwall's inequality, we obtain
$$(\e+\|u_t\|_{\C_b^1})^{\theta^{-1}} \le c_2\e^{c_2\int_0^t \log(\e+ \|u_s\|_{\C_b^1})   \d r},\ \ t\in [0, T_n], n\ge 1.$$
Thus, there exists a constant $c_3>0$ such that
$$\log (\e+ \|u_t\|_{\C_b^1})\le c_3+ c_3\int_0^t \log(\e+ \|u_s\|_{\C_b^1})   \d r,\ \ t\in [0,T_n], n\ge 1.$$
Using Gronwall's inequality again, we derive
$$\log (\e+ \|u_t\|_{\C_b^1})\le c_3\e^{c_3  T},\ \ t\in [0,  T_n], n\ge 1.$$
Therefore, there exists a constant $c_4>0$ such that
\beq\label{AD} \sup_{t\in [0,  T_n]}  \|u_t\|_{\C_b^1} \le c_4,\ \ n\ge 1.\end{equation} 
By \eqref{TN},  if $T_n<T$ then this implies 
$n\le c_4$. Thus, $T_n\ge T$ for $n> c_4$. Hence  $T^*=T_n=T$  for large $n\ge 1$, so that \eqref{NN} follows from \eqref{*PP0} and \eqref{AD}. 

\end{proof}



{\footnotesize
\begin{thebibliography}{999}






\bibitem{BR} V. Barbu, M. R\"ockner, \emph{From nonlinear Fokker-Planck equations to solutions of distribution dependent SDE,}   Ann. Probab. 48(2020), 1902--1920.

\bibitem{Bismut} J. M. Bismut, \emph{Conjugate convex functions in optimal stochastic control,} J. Math. Anal. App. 44(1973), 384--404. 

\bibitem{CSTV} P. Cheridito, H. M. Soner, N. Touzi, N. Victoir, \emph{Second-order backward stochastic differential equations and fully nonlinear parabolic PDEs,} 
Comm. Pure Appl. Math. 60(2007), 1081--1110.

\bibitem{CI} P. Constantin, G. Iyer, \emph{A stochastic representation of the three-dimensional incompressible Navier-Stokes equations,} Comm. Pure Appl. Math. 61(2008), 330--345. 

\bibitem{FF} U. Frisch, \emph{Turbulence. The legacy of A. N. Kolmogorov,}  Cambridge University Press, Cambridge (1995).

\bibitem{Gi}   B. Gilding, M. Guedda, R. Kersner, \emph{The Cauchy problem for $  u_t=\DD u + |\nn u|^q$,}
 J. Math. Anal. Appl. 284(2003), 733--755.

 
 \bibitem{Ito} K. It\^o, \emph{Differential equations determining a Markoff process,} J. Pan-Jpn. Math. Coll. No. 1077(1942). In Kiyosi It\^o: Selected Papers, Springer (1987). 
 

\bibitem{KL} A. N. Kolmogorov, \emph{Local structure of turbulence in an incompressible fluid at very high Reynolds numbers,} Dokl. Akad. Nauk SSSR 30(1941), 299--303.

\bibitem{KS}  J. Krug, H. Spohn, \emph{Universality classes for deterministic surface growth,} Phys. Rev. A(3)
38(1988), 4271--4283.

\bibitem{MC} H. P.  McKean,  \emph{
   A class of Markov processes associated with nonlinear parabolic equations, }  Proc. Nat. Acad. Sci. U.S.A. 56 (1966), 1907--1911.

\bibitem{PP} E. Pardoux, S. Peng, \emph{Adapted solutions of a backward stochastic differential equation,} Syst. Control Lett. 14(1990), 55-61.

\bibitem{Peng0} S. Peng, \emph{$G$-expectation, $G$-Brownian motion and related stochastic calculus of It\^o's type,} in: Benth et al. (eds.) The Abel Symposia, pp. 541--567. Springer (2007). 
\bibitem{Peng} S. Peng, \emph{Nonlinear Expectations and Stochastic Calculus under Uncertainty,} Springer 2019. 

\bibitem{Qian}   Z. Qian,  E.  S\"uli, Y.  Zhang,   \emph{Random vortex dynamics via functional stochastic
differential equations,}  arXiv:2201.00448v1.
\bibitem{W21e} F.-Y., Wang,  \emph{Regularity estimates and intrinsic-Lions derivative formula   for  singular McKean-Vlasov  SDEs,}   arXiv:arXiv:2109.02030.
\bibitem{22} F.-Y. Wang, \emph{Distribution dependent SDEs  for   Navier-Stokes type equations,}   Electron. Commun. Probab. 27(2022),   1--12.
 \bibitem{XXZZ}  P. Xia,  L.   Xie, X.  Zhang, G. Zhao, \emph{$L^q$($L^p$)-theory of stochastic differential equations,}   Stoch. Proc. Appl. 130(2020), 5188--5211.


\bibitem{YZ}   C.  Yuan, S.-Q. Zhang, \emph{A study on Zvonkin's transformation for stochastic differential equations with singular drift and related applications,} J. Diff. Equat.   
297(2021),  277--319.

 \end{thebibliography}}
\end{document}